\documentclass[times,final]{article}

\usepackage{amssymb}
\usepackage{latexsym}
\usepackage{array,amsfonts,amsmath,amssymb,graphicx,relsize,url,array,indentfirst,booktabs,geometry,subfig}
\usepackage{listings}
\usepackage{amsthm}
\usepackage{titling}
\usepackage{mathtools,color}
\usepackage{url}
\usepackage{xcolor}
\definecolor{newcolor}{rgb}{.8,.349,.1}

\begin{document}

		\title{A novel multi-scale loss function for classification problems in machine learning}
		
		\author{Leonid Berlyand \thanks{Department of Mathematics, The Pennsylvania State University, University Park, Pennsylvania, USA}
			\and
		Robert Creese\protect\footnotemark[1]
		\and
		Pierre-Emmanuel Jabin\protect\footnotemark[1]
	 }

		\maketitle
		

		\begin{abstract}
			We introduce two-scale loss functions for use in various gradient descent algorithms applied to classification problems via deep neural networks.  This new method is generic in the sense that it can be applied to a wide range of machine learning architectures, from deep neural networks to support vector machines for example. These two-scale loss functions allow to focus the training onto objects in the training set which are not well classified.  This leads to an increase in several measures of performance for appropriately-defined two-scale loss functions with respect to the more classical cross-entropy when tested on traditional deep neural networks on the MNIST, CIFAR10, and CIFAR100 data-sets.
		\end{abstract}
		

	
	
	\section{Introduction}
	We introduce in this article a novel multi-scale perspective on loss functions that are critical for the training phase of machine learning algorithms. The main idea of those additional scales is to allow to focus the training on the objects in the datasets that are still poorly handled by the algorithms. As a first step in the present article, we focus on two-scale loss functions.  One advantage of our multi-scale loss functions is that they can generically be applied to wide range of problems and techniques for learning algorithms, as they do not require any change in the structure of the method but only a straightforward modification of how the loss or cost is estimated. 
	However, for simplicity in this article, we describe and test the methods on classical classification problems based on deep neural networks techniques.
	
	Deep neural networks (DNNs) are commonly used for classification problems, given their demonstrated performance (see for example the review \cite{rawwan17}).  To briefly summarize, the goal of a classification algorithm is to predict the class, $i(s)$, of each object in $s$ in a given dataset.  Deep learning algorithms typically lead to statistical classifiers with an output of predicted probabilities, $p=(p_1(s), p_2(s), \dots, p_n(s))$, where $p_i(s)$ is the predicted probability that $p$ belongs to the $i$-th class.  One then usually assigns the index of the largest probability as the predicted class.
	
	To evaluate the performance of a classifier, one makes use of a so-called test set $T\in S$, on which the correct class is know. There are various measures of performance, which may have different behaviors. The overall accuracy of $T$,
	\[
	acc_T(\phi(\alpha))=\dfrac{\#\{s\in T: p_{i(s)}=\max_{1\leq i\leq K} p_i\}}{\#T},
	\]
	is a common choice, together with other measures such as the average accuracy which estimates the accuracy on each class and then takes the average.		
	Overall accuracy (together with other accuracy measures) has the drawback of not considering the actual predicted probabilities and for this reason does not reveal the full picture of DNN performance.   For example, consider two objects $s_1$ and $s_2$ and suppose there are only 2 classes. If $p_{i(s_1)}=0.51$ and $p_{i(s_2)}=0.99$, the object $s_2$ is clearly classified better, but both would just be classified as correct and have the same impact on overall accuracy.  Alternatively, if $p_{i(s_1)}=0.01$ and $p_{i(s_2)}=0.49$, both objects are classified incorrectly and neither object would contribute to the accuracy.  However the classifier is clearly doing a better job on $s_2$.  From this example, we see the need for other measures of performance which can take into account how close objects are to being correctly classified with respect to probability.
	
	One such example is  top-$k$ performance, which counts the objects whose probabilities for the correct class are one of the $k$ largest. Top-$k$ is used  frequently in image classification, e.g. the CIFAR100 data-set and Imagenet~\cite{KriHin09, KriSutHin12}.  In addition, the CIFAR100 data-set in particular has several classes which are similar, so an alternative set of ``super"-classes, groups of 5 similar classes, are used to measure performance as well \cite{KriHin09}. 
	We will also make use of another measure that estimates ``close-enough'' performance to better capture those objects that are close to being well classified.

	When attempting to improve the performance of a classifier, one can try to adjust the structure, e.g. number, width, and type of layers, of the DNN classifier  according to the specific classification problem. 
	Several DNNs have already incorporated  multi-scale structures in order to improve performance of certain problems, e.g. in image classification \cite{WanDuZha20}, spacecraft control \cite{CheWanJia19},  biology \cite{SamVee19}, and the MscaleDNN network in PDEs \cite{LiXUZha20}.  In all of these cases the structure of the DNN, described in Sec. \ref{Sec:2_scale_loss}, is modified in order to extract multi-scale features which appear naturally in these problems.  For example, in image classification, features are found at different scales, corresponding for instance to the contours of objects of very different sizes. 
	These multi-scale neural networks aim to ensure that all of these features are each captured by creating the appropriate scales in the structure of their network.  Additionally multi-scale average loss functions are studied, where the scales appear through variations in the data set and due to the network structure \cite{KonTao20}.
	
	On the other hand, techniques not specific to a type of  classification problem are used as well, such as developing new training procedures, e.g. ADAM optimization for DNNs \cite{KinBa14}, or choosing alternative loss functions, as in \cite{deBPas15}.  
	Our multi-scale loss function is in that spirit 
	as it leaves the structure of the network fully unchanged.
	That makes it applicable to a wide variety of DNNs used for classification, including those with multi-scale structures. 
	
	We note that in general the connection between changes in the loss function and the performance on the test set is delicate and not immediate. {A first issue results from the fact that the performance can be measured both on the training set $TR$ and the test set $T$, which are typically separate, $T\cap TR=\emptyset$.}  
	In general increasing performance on the training set does lead to an increase of accuracy on the test set in practice, however that outcome is not guaranteed.  Even on the training set, there may not always be an immediate connection between a decrease in the loss function leads and an increase in accuracy, leading to various questions of stability for study as in \cite{BerJabSaf20}.

	{We describe a basic structure of the type of DNNs that we use in Sec. \ref{Sec:2_scale_loss}}. We present in the same section our two-scale loss function and explain the procedure to select its scales in order to focus training onto objects which are close to being correctly or incorrectly classified (e.g. $p_{i(s)}\approx 0.5$ if there are 2 classes).  \textcolor{black}{In Sec. \ref{Sec:Num} we compare the two-scale loss function against a traditional single scale loss function for various measures of performance for several data-sets, MNIST, CIFAR10, and CIFAR100.  
	We see large and clear improvements from a small increase of computational complexity for each of the data-sets tested for appropriate measures of performance. }

	\section{Two-scale loss functions}
	\label{Sec:2_scale_loss}
	\subsection{Review of DNN Notation}
	We first give a brief description of the classification problem. A set of objects $S\in \mathbb{R}^n$, is considered.  Each object $s\in S$ belongs to exactly one of $K$ classes where each class is denoted by an integer between $1$ and $K$.  The class to which an object belongs is called the correct class and denoted by an integer $1\leq i(s) \leq K$.  However the correct class is only known on a finite subset $TR$ of $S$, called the training set.  The aim of the classification problem is to approximate $i(s)$ on all objects in $S$, by using the training set where $i(s)$ is known.  To approximate $i(s)$, an family of classifiers $\phi(\alpha,s)$, defined by the parameter $\alpha\in \mathbb{R}^\mu$ is considered. This family of classifiers has an output of a finite probability distribution,
	\[		\phi(\alpha,s)=(p_1(\alpha,s),p_2(\alpha,s),\dots,p_K(\alpha,s)),\quad \sum_{i=1}^K p_i=1,\quad 0\leq p_i\leq 1,
	\]
	where each $p_i$ denotes the probability of belonging to the $i$th class.  This is a ``soft" or statistical classifier which gives a probability distribution as an output instead of a class.  We can have the DNN predict which class an object belongs to by taking 
	\[
	\arg \max_{1\leq i\leq K} p_i,
	\]
	as the predicted class.  Since $i(s)$ denotes the correct class, if the predicted class is $i(s)$ then we say that the object is correctly classified by the DNN.

	To describe our original motivation for a two-scale loss function, some basic notations for deep neural networks 
	for the classification problem are introduced. 
	%
	We consider a DNN as a composition of $M$ layers, where the first $M-1$ layers are compositions of a linear operation followed by a non-linear operation,  \begin{align} Y^{l+1}_i=\sum\limits_{j=1}^{N_{l}} \alpha_{i,j}^{l+1} X^{l}_{j},\quad l=0\dots M-2, \quad {{i = 1, 2, \cdots N_{l + 1}}}\label{connectedl} \\
		X^{l+1}_i=\lambda\left(Y^{l+1}_i\right),\quad l=0\dots M-2, \quad {{i = 1, 2, \cdots N_{l + 1}}}.\label{lambdal}.
	\end{align} 
	The value $N_l$ denotes the dimension of the $l$th layer and $X^{l}\in \mathbb{R}^{N_l}$ is the vector of values at the $l$th layer.  The terms $\alpha_{i,j}^{l+1}$ are defined by the parameter $\alpha$ and define the linear operation.  The function $\lambda:\mathbb{R}\to\mathbb{R}$ is a nonlinear activation function, with a popular choice being the Rectified Linear Unit (ReLU), $\lambda(x)=\max(x,0)$ \cite{LecBenHin15}.  For the final layer we consider a probability normalization function, $h:\mathbb{R}^K\to\mathbb{R}^K$, satisfying
	\[
	\sum\limits_{i=1}^K h\left(x\right)_i=1 \text{ and } 0\leq h\left(x\right)_i\leq 1 \, \forall 1\leq i\leq K.
	\]
	and an order-preservation property 
	\[
	x_i<x_j \implies h\left(x\right)_i<h\left(x\right)_j\, \text{ for all } 1\leq i,j\leq K.    
	\]
	While our motivation applies to any probability normalization function with order preservation, we  will focus on the frequently used soft-max function,
	\begin{equation}
		X^M_i:=p_i\left(\alpha,s\right)=\dfrac{e^{X^{M-1}_i}}{\sum\limits_{j=1}^K e^{X^{M-1}_j}},\, i=1\dots K,
		\label{softmax}
	\end{equation}
	in order to clearly explain our motivation for the two-scale loss function. 
	
	The structure of a DNN refers to choice of the number of layers, the dimension of each layer, and the type of linear map chosen.  For example the matrix $a^{l+1}=[a_{ij}^{l+1}] \in \mathbb{R}^{N_l \times N_{l+1}}$ defining the linear layer may have restrictions, e.g. $a^{l+1}$ can be restricted to be a Toeplitz matrix, which creates a type of layer known as a convolutional layer \cite{jain89}.  If there are no restrictions on $a^{l+1}$ then the layer is called fully connected.  Additionally a layer may have no parameters defining \eqref{connectedl}, and may instead have a predefined function, e.g. taking the maximum of every 4 elements (assuming $N_l$ is divisible by 4),
	\[
	Y^{l+1}_i=\max\{x_{4i+1},x_{4i+2},x_{4i+3},x_{4i+4}\},\quad i=1,2,\cdots,N_l/4=N_{l+1},
	\]
	which is an example of a layer known as a max pooling layer.  Much more complex DNN structures have been introduced, but we will focus on DNNs with layers defined via \eqref{connectedl}-\eqref{lambdal} for the motivation of our two-scale network.
	
	To choose a good classifier from the family $\phi(\alpha,s)$, {\em i.e.} to identify a good choice of parameter $\alpha$, one introduces a so-called loss function defined on each object;   a typical choice being the cross entropy,
	\begin{equation}
		L(\alpha,s)=-\log (p_{i(s)}(\alpha,s)). \label{cross_entropy}
	\end{equation}
	The average loss over the training set can then be defined by
	\begin{equation}
		\bar L(\alpha)=\dfrac{1}{\#TR}\sum_{s\in TR} L(\alpha,s).\label{trad_loss}
	\end{equation}
	Note that the minimum $\bar L(\alpha)=0$ is reached if and only if $p_{i(s)}=1$ for each object $s$, that is if the classifier is perfect.
	
	The identification of an appropriate set of parameter $\alpha$ is performed during the so-called training phase. The goal of the process is find $\alpha$ with $\bar L(\alpha)$ as small as possible. Various optimization algorithms can be used in that regard; early on, iterative gradient descent (GD) algorithm,
	\begin{equation}
		\alpha^{n+1}=\alpha^{n}-\tau\nabla_\alpha \bar L\left(\alpha^n\right), \label{GD}
	\end{equation}
	were commonly employed.  Here $n$ denotes the step of training and $\alpha^0$ is the random initial choice of parameters.
	
	However, as the training set can be rather large (and is ideally very large to ensure good performance), this leads to a large computational cost for each step of \eqref{GD}.  For this reason, variants of GD such as the stochastic gradient descent (SGD) are usually preferred nowadays. SGD for example can be seen as using a randomly sampled loss, instead of the full averaged loss $\bar L(\alpha,s)$, which is given at each step $n$ by
	\[
	\hat L_n=  L(\alpha,s)=\dfrac{1}{\#B}\sum_{s\in B_n}-\log (p_i(s)),
	\]
	where $B_n$ is set of $\#B$ randomly selected objects from the training set.  The term $\#B\in \mathbb{N}$ is called the batch size. How to choose the initial parameter value $\alpha_0$ can be another delicate point; random initial $\alpha_0$ have proven popular with for example each term in $\alpha_0$ chosen from a Gaussian distribution. 
	\subsection{Motivation for a two-scale loss function}
	We recall that, if the output of the DNN for an object $s$ satisfies,
	\[
	p_{i\left(s\right)}\left(\alpha,s\right)>\max_{1\leq i\leq K, i\neq i\left(s\right)} p_i\left(\alpha,s\right),
	\]
	then we refer to $s$ as correctly classified since the probability of belonging to the correct class $i\left(s\right)$ is the highest. Additionally since soft-max is order preserving, $s$ is correctly classified if \begin{equation}
		X_{i\left(s\right)}^{M-1}\left(\alpha,s\right)>\max_{1\leq i\leq K, i\neq i\left(s\right)} X_{i}^{M-1}\left(\alpha,s\right).\label{pen_ult_class}
	\end{equation}
	We also consider that the linear activation function is often homogeneous of degree one, which is the case for the ReLU or absolute value function. This results from each layer \eqref{connectedl}-\eqref{lambdal} being homogeneous of degree one with respect to $\alpha$ as a composition of a linear function and a homogeneous activation function. Each $X_{i}^{M-1}\left(\alpha,s\right)$ is then homogeneous of degree $M-1$ with respect to $\alpha$,
	\begin{equation}
		X_{i}^{M-1}\left(c\alpha,s\right)=|c|^{M-1}X_{i}^{M-1}\left(\alpha,s\right).   \label{homogenous}
	\end{equation}
	Observe that if we multiply the parameter $\alpha$ in \eqref{pen_ult_class} by some scalar $c$, we can use the homogeneity to show that 
	\[
	X_{i\left(s\right)}^{M-1}\left(c\alpha,s\right)>\max_{1\leq i\leq K, i\neq i\left(s\right)} X_{i}^{M-1}\left(c\alpha,s\right)\Longleftrightarrow |c|^{M-1}X_{i\left(s\right)}^{M-1}\left(\alpha,s\right)>|c|^{M-1}\max_{1\leq i\leq K, i\neq i\left(s\right)} X_{i}^{M-1}\left(\alpha,s\right),
	\]
	which is equivalent to \eqref{pen_ult_class} without the scaling term $c$ for the parameter $\alpha$.  As a result which class has the largest probability, second largest probability, etc., does not depend on the size of the parameters $\alpha$, which we may simply scale as 
	\begin{equation}
		R=\left(\Pi_{i=1}^{M-1} \| (\alpha^i)\|\right),\label{R_def}
	\end{equation}
	where $\alpha^i$ is the vector of all parameters defining the linear layer $\#i$ in \eqref{connectedl}. The norm $\|\cdot\|$ is the Euclidean norm.  Note that $R$ is positive so we do not need the absolute value sign, when applying the homogeneity property to $R$.  
	
	Since $R$ does not affect the order of the probabilities, scaling the parameters does not affect the accuracy, top-$k$ accuracy, or Rand index for example.  These measures of performance
	only depend on the direction of $\alpha$ which we define as follows
	\[
	\hat \alpha=\frac{\alpha}{R^{\frac{1}{M-1}}}.
	\]	
	This leads us to introduce the normalized values at the layer before soft-max by
	\[
	\hat X_i^{M-1}\left(\alpha\left(t\right),s\right):=X_i^{M-1}\left(\hat\alpha\left(t\right),s\right).
	\]
	To see how $R$ influences the loss function we rewrite \eqref{cross_entropy} in terms of $\hat X^m$ in order to see the influence of the size of the parameters, $R$.  First note that by combining the soft-max function \eqref{softmax} with the cross-entropy loss function \eqref{cross_entropy}, one obtains 
	\[
	L\left(\alpha,s\right)=-\log\left(\frac{e^{\left( X_{i\left(s\right)}^{M-1}\left( \alpha,s\right)\right)}}{\sum\limits_{1\leq i\leq K}e^{\left( X_{i}^{M-1}\left( \alpha,s\right)\right)}}\right)=-\log\left(\frac{e^{\left(R \hat X_{i\left(s\right)}^{M-1}\left( \alpha,s\right)\right)}}{\sum\limits_{1\leq i\leq K}e^{\left(R\hat X_{i}^{M-1}\left( \alpha,s\right)\right)}}\right). 
	\]
	Dividing both the numerator and denominator inside the logarithm by the numerator yields
	\begin{equation}
		L\left(\alpha,s\right)=-\log\left(\frac{1}{1+\sum\limits_{1\leq i\leq K, i\neq i\left(s\right)}e^{R\left(\hat X_{i}^{M-1}\left( \alpha,s\right)-\hat X_{i\left(s\right)}^{M-1}\left( \alpha,s\right)\right)}}\right). \label{loss_dX}
	\end{equation}
	If an object is classified correctly, then $X_{i\left(s\right)}^{M-1}\left( \alpha,s\right)$ is the largest component of $X^{M-1}\left( \alpha,s\right)$.  Then each exponent in \eqref{loss_dX} is negative, showing that if $R\to \infty$ then $L\left(\alpha,s\right) \to 0$ and all correctly classified objects are perfectly classified.  However sending $R\to \infty$ would make all objects which are not correctly classified, even those with classification confidence near zero, very poorly classified.  If all objects are correctly classified after some time then it is possible to show that $R\to \infty$, see \cite{BerJab18}.
	
	Next we define the classification confidence of an object in the test set $T$ or training set $TR$ as
	\[
	\delta X\left( \alpha,s\right):=X_{i\left(s\right)}^{M-1}\left( \alpha,s\right)-\max_{i\neq i\left(s\right)}X_{i\left(s\right)+1}^{M-1}\left( \alpha,s\right). 
	\]
	and the normalized classification confidence as
	\[
	\delta \hat X\left( \alpha,s\right):=X_{i\left(s\right)}^{M-1}\left( \hat \alpha,s\right)-\max_{i\neq i\left(s\right)}X_{i\left(s\right)+1}^{M-1}\left( \hat \alpha,s\right)=\delta X\left(\hat{\alpha},s\right)=\dfrac{\delta X\left(\alpha,s\right)}{R}. 
	\]
	If $\delta X\left( \alpha,s\right)>0$ then $X_{i\left(s\right)}^{M-1}\left( \alpha,s\right)$ would be the largest component of $X^{M-1}\left( \alpha,s\right)$. Therefore the probability of belonging to the correct class, $p_{i\left(s\right)}\left(\alpha,s\right)$, would be the largest and $s$ would be correctly classified. Similarly if $\delta X\left( \alpha,s\right)<0$, then the maximum probability is not $p_{i\left(s\right)}\left(\alpha,s\right)$ and $s$ is incorrectly classified.  Additionally $\delta X\left( \alpha,s\right)$ is uniformly Lipschitz continuous with respect to $s$ as a composition of Lipschitz functions \eqref{connectedl}-\eqref{lambdal}. Therefore a small change in $s$ will not change the classification of $s$. 
	
	We now introduce the set of well classified objects with confidence $\eta\geq 0$
	\[
	W_{\eta,I}\left(\alpha\right)={s \in I: \;  \delta X\left( \alpha,s\right)>\eta}, \quad I\subset S.
	\]
	The accuracy of the DNN with parameter $\alpha$ can be recovered immediately with
	\begin{equation}
		acc\left(\alpha\right)=\dfrac{\#W_{0,T}\left(\alpha\right)}{\#T}.    \label{acc_defn}
	\end{equation}
	The dependence of $\#W_{\eta,T}$ on $\eta$ and more generally the distribution of $\delta X(\alpha,s)$ with respect to $s$ provides a significantly larger amount of information about the performance of the network than other measures of performance, e.g. overall accuracy. For example if $\#W_{\eta,T}$ is large for $\eta=0$ but sharply decreases as $\eta$ grows, we can infer that many objects are just barely classified correctly.

	It has already been recognized that cross-entropy loss functions  may not always be the best and there have been many proposed changes to the cross-entropy loss function in an attempt to improve accuracy or resistance to noise in the training set; see \cite{ZhaSab18,MarSti18,PanSunKha16}.  For example, the paper~\cite{ZhaSab18} introduced a truncated loss function as
	\begin{equation}
		L_{Trunc}\left( \alpha,s\right) =\begin{cases} -\log(k) &\mbox{if }  p_{i(s)}<k \\
			L\left( \alpha,s\right) &\mbox{if }  \delta X\left( \alpha,s\right)\geq k \end{cases}, \label{trunc_loss}
	\end{equation}
	for $0 \leq k\leq 1$.  If an object has a very low probability for the correct class the loss for that single object cannot exceed $-\log(k)$.  As a result a single object can not contribute a very large amount of loss, thereby mitigating the effect an outlier or mislabeled object might have on training.  If $k=0$, $L_{Trunc}$ reverts to the traditional loss function as $p_{i(s)}$ can not be less than $0$.  This choice of loss function increased the performance of the chosen DNN on several data-sets where noise was added, though accuracy was somewhat lower when noise was not added.  

	The main idea in this article is to instead  focus training onto objects which are not yet well classified. For objects in $W_{\eta,TR}\left(\alpha\right)$, if $R$ is large, the exponent for each term in \eqref{loss_dX} decreases, and pushes $L\left(\alpha,s\right)$ closer to $0$ irrespective of $\hat\alpha$.  Therefore if $R$ is increased only for objects $\delta X(\alpha,s)>\eta$, the loss decreases on the well classified objects and training should focus on the other, not well classified objects.
	
	To this end,
	the introduced loss functions will split the training set into two cases, one for the set of well classified objects in the training sets, $W_{\eta,TR}\left(\alpha\right)$, and another for the remaining objects.

	
	\subsection{Two-scale loss functions}\label{Sec:2_loss_var}
	This leads us to the definition of the two-scale loss function as
	\begin{equation}
		\label{2_scale_loss} L_{2,\eta}\left(\left( \alpha,R_1,R_2\right),s\right) = \begin{cases} -\log\left(\dfrac{1}{1+\sum\limits_{1\leq i\leq K, i\neq i\left(s\right)}e^{R_1\left(\hat X_{i}^{M-1}\left( \alpha,s\right)-\hat X_{i\left(s\right)}^{M-1}\left( \alpha,s\right)\right)}}\right) &\mbox{if }  \delta X\left( \alpha,s\right)<\eta \\
			-\log\left(\dfrac{1}{1+\sum\limits_{1\leq i\leq K, i\neq i\left(s\right)}e^{R_2\left(\hat X_{i}^{M-1}\left( \alpha,s\right)-\hat X_{i\left(s\right)}^{M-1}\left( \alpha,s\right)\right)}}\right) &\mbox{if }  \delta X\left( \alpha,s\right)\geq \eta \end{cases}, 
	\end{equation}
	where $\eta>0$ is a chosen, fixed parameter. $R_1$ and $R_2$ are two new parameters which denote the two scales of the parameters used in \eqref{2_scale_loss} and can be adjusted during training.  With new parameters to train, the deterministic gradient descent algorithm becomes
	\begin{equation}
		\left( \alpha^{n+1},R_1^{n+1},R_2^{n+1}\right)=\left( \alpha^n,R_1^n,R_2^n\right)-\tau \nabla_{\left( \alpha,R_1,R_2\right)} \sum\limits_{s\in TR} L_{2,\eta}\left(\left( \alpha^n,R_1^n,R_2^n\right),s\right), \label{GD_2s}
	\end{equation}
	with corresponding formula for stochastic gradient algorithms.
	
	The desired quality of decreasing loss for well classified objects is achieved with \eqref{2_scale_loss} when $R_2\gg R_1$ and $\eta>0$.  We note that, if we chose $\eta<0$, we cannot predict the effect of a large $R_2$ on the loss in general.  However, for example, if there are only 2 classes, then $R_2\gg R_1$ would increase the loss for the objects with $\eta<\delta X(\alpha,s)<0$.  As a result these objects would be the focus of more training. This contradicts the original goal of the two-scale loss function, focusing training on the objects with lower classification confidence, $\delta X(\alpha,s)<\eta$.  Therefore we restrict ourselves to considering $\eta>0$.
	
	\textcolor{black}{We chose to have large $R_2$ as for well classified objects the tendency is to push $R$ to grow throughout training.  Instead of focusing  on making the size of the parameters larger, we aim to have the training procedure focus on improving the direction of the parameter by choosing a larger scale for the well-classified objects.}
	
	While we perform our testing in the next section for the loss function defined by \eqref{2_scale_loss}, it can be useful to observe that any loss function into a two-scale loss function using the form
	\[
	L_{2,\eta}\left(\left( \alpha,R_1,R_2\right),s\right) =\begin{cases} L\left( (\hat \alpha,R_1),s\right) &\mbox{if }  \delta X\left( \alpha,s\right)<\eta \\
		L\left( (\hat\alpha,R_2),s\right) &\mbox{if }  \delta X\left( \alpha,s\right)\geq \eta \end{cases},
	\]
	where $R_1$ and $R_2$ are used to replace the size of the parameters. This will similarly decrease loss on the well classified objects when $R_2>R_1$ if the chosen loss function decreases when the probability of the correct class increase and all other probabilities either decrease or increase less than the correct probability. For example the original purpose of decreasing loss on well-classified objects works on the mean squared error loss as well as alternative loss functions used in \cite{ZhaSab18,MarSti18} as well as alternatives to the soft-max in \cite{marAst16,deBPas15,AsaLit17}, e.g. sparsemax.
	

	As claimed in the introduction, the two-scale loss function only depends on the loss function and normalizing function, e.g. soft-max.  As a consequence, no matter the structure of the network that yields the penultimate layer $X^{M-1}$, one can use this two-scale loss function \eqref{2_scale_loss}.
	There are consequently a wide variety of DNN structures, data-sets, and loss functions where this novel modification can be applied, which we plan to explore in future works.  
	
	
	We can infer some elementary behavior of $R_2$ during training by considering two cases.  The first case is where no objects in the training set $TR$ are well-classified.  In that case $R_2$ is never used in calculating the loss, \eqref{2_scale_loss}, so the partial derivative $\frac{\partial}{\partial R_2}L_{2,\eta}\left(\left(\hat \alpha,R_1,R_2\right),s\right)$ is $0$ for all $s \in TR$ and $R_2$ does not change.  The other case is when at least one object is well classified.  If an object $s$ is well classified the partial derivative of loss with respect to $R_2$ for that object is given by
	\begin{equation}
		\dfrac{\partial}{\partial R_2}L_{2,\eta}\left(\left(\hat \alpha,R_1,R_2\right),s\right)=\frac{\sum\limits_{1\leq i\leq K, i\neq i\left(s\right)}\left(\hat X_{i}^{M-1}\left( \alpha,s\right)-\hat X_{i\left(s\right)}^{M-1}\left( \alpha,s\right)\right)e^{R_2\left(\hat X_{i}^{M-1}\left( \alpha,s\right)-\hat X_{i\left(s\right)}^{M-1}\left( \alpha,s\right)\right)}}{1+\sum\limits_{1\leq i\leq K, i\neq i\left(s\right)}e^{R_2\left(\hat X_{i}^{M-1}\left( \alpha,s\right)-\hat X_{i\left(s\right)}^{M-1}\left( \alpha,s\right)\right)}}. \label{R2_der}
	\end{equation}
	The classification confidence of the object is greater than zero so by 
	\[
	\hat X_{i}^{M-1}\left( \alpha,s\right)-\hat X_{i\left(s\right)}^{M-1}\left( \alpha,s\right)\leq\max_{1\leq i\leq K, i\neq i\left(s\right)}\hat X_{i}^{M-1}\left( \alpha,s\right)-\hat X_{i\left(s\right)}^{M-1}\left( \alpha,s\right)
	=-\delta \hat X\left(\alpha,s\right)<0,
	\]
	the term $\hat X_{i}^{M-1}\left( \alpha,s\right)-\hat X_{i\left(s\right)}^{M-1}\left( \alpha,s\right)$ must be negative.   Therefore the expression \eqref{R2_der} is negative for every well classified object $s$.  If $s$ is not well classified the partial derivative is $0$ as $R_2$ is not calculated in loss.  This results in $R_2$ growing  at each step, since the parameters move in the direction of the negative gradient during training.  Therefore the DNN will have $R_2$ grow throughout training, as long as one object is well classified.  This is reminiscent of a result in \cite{BerJab18} where if all objects are well classified the magnitude of the parameters $\alpha$ will grow.  However with the two scale training given by \eqref{GD_2s}, we are unable to specify the rate of growth of $R_2$, contrary to \cite{BerJab18}.
	

	One should also remark that the two-scale loss function \eqref{2_scale_loss} may return to a one-scale case if $R_2=R_1$ after some training, since $R_1$ may grow faster than $R_2$ for example.  Alternatively if every object in $TR$ becomes well classified or every object in $TR$ becomes miss-classified, then only one of $R_1$ or $R_2$ is used and a one-scale loss function is also recovered for all practical purposes. This should ensure that training with the loss function \eqref{2_scale_loss} will perform at least as well as the traditional loss function \eqref{cross_entropy} as long as appropriate choices are made for $\eta$ and initial $R_1$ and $R_2$.  
	
	We will also consider a variant of the two-scale loss function where the scales $R_1$ and $R_2$ are fixed.  The loss function is still given by \eqref{2_scale_loss}, however $R_1$ and $R_2$ are no longer parameters and instead chosen scalars.  The gradient descent training algorithm would then change back from \eqref{GD_2s} to
	\begin{equation}
		\alpha^{n+1}= \alpha^n-\tau \nabla_{ \alpha} \sum\limits_{s\in TR} L_{2,\eta}\left(\left( \alpha,R_1,R_2\right),s\right). \label{GD_fixed_2s}
	\end{equation}
	The above method of training with \eqref{2_scale_loss} and \eqref{GD_fixed_2s} will be referred to as fixed two-scale loss for convenience.  The desired quality of decreasing loss for well classified objects is again achieved with \eqref{GD_fixed_2s} when $R_2\gg R_1$ and $\eta>0$. 
	
	An advantage the fixed two-scale loss has is that the scales $R_1$ and $R_2$ can be ensured to remain well separated, as they maintain their initial values for all time.  
	The partial derivatives with respect to $R_1$ and  $R_2$ no longer need to be calculated, leading to a small decrease in computational complexity.
	
	But a major concern with the fixed two-scale loss function however is that training can be slowed due to the fixed values of $R_1$ and $R_2$, as we will indeed observe in numerical testing.  It is often advantageous to vary the size of the parameters during training, and  the fixed two-scale loss function is unable to capture that advantage while the algorithm in \ref{Sec:2_loss_var} does.  Additionally a poor choice of values for $R_1$ and $R_2$ is disadvantageous and cannot be corrected via training with the fixed two-scale loss function, whereas the algorithm using \eqref{GD_2s} would improve upon the choice of values throughout time. 
	
	Having both scales fixed is hence restrictive and for this reason, we also introduce a third variant where $R_1$ remains fixed but the ratio $R_2/R_1$ is allowed to change, leading to the following
	\begin{equation}
		\label{2_scale_loss_sep}
		L_{2_{sep},\eta}\left(\left( \alpha,R_s\right),s\right) = \begin{cases} -\log\left(\dfrac{1}{1+\sum\limits_{1\leq i\leq K, i\neq i\left(s\right)}e^{R\left(\hat X_{i}^{M-1}\left( \alpha,s\right)-\hat X_{i\left(s\right)}^{M-1}\left( \alpha,s\right)\right)}}\right) &\mbox{if }  \delta X\left( \alpha,s\right)<\eta \\
			-\log\left(\dfrac{1}{1+\sum\limits_{1\leq i\leq K, i\neq i\left(s\right)}e^{R_s R\left(\hat X_{i}^{M-1}\left( \alpha,s\right)-\hat X_{i\left(s\right)}^{M-1}\left( \alpha,s\right)\right)}}\right) &\mbox{if }  \delta X\left( \alpha,s\right)\geq \eta \end{cases} 
	\end{equation}
	where $R$ is fixed but the term $R_s$ denotes how well separated the scales are and evolves with training.  Correspondingly, the gradient descent must be modified to 
	\begin{equation}
	\left(\alpha^{n+1},R_s^{n+1}\right)=\left( \alpha^n,R_s^n\right)-\tau \nabla_{\left( \alpha,R_s\right)} \sum\limits_{s\in TR} L_{2_{sep},\eta}\left(\left( \alpha^n,R_s^n\right),s\right), \label{GD_2s_sep}
	\end{equation}
	since $R_s$ is a trainable parameter.  
	
	Notice that the partial derivative of the loss with respect to $R_s$ is strictly non-negative by the same reasoning that $R_2$ is increasing in \ref{Sec:2_loss_var}. Therefore if $R_s$ evolves via training, the scales can only become ever more separated.  This is not necessarily the case in our initial loss function \ref{Sec:2_loss_var}: The growth of $R_1$ may be larger than the growth of $R_2$, leading $R_2/R_1$ to shrink.  This third loss function is more restrictive but can still work well when a separation of scales is desirable.
        \subsection{Theoretical conjectures}
        \textcolor{black}{ 	While in this paper we do not address theoretical justification of the 2-scale loss function introduces in Sec. \ref{2_scale_loss}, we nevertheless formulate several conjectures for our future works.  
        	First, given a network and data sets (at least a training set), one would like to find the minimizant $\alpha^*$ for the classical loss function $L(\alpha)$ together with the minimizant $\alpha^*_\eta$ for the 2-scale loss function $L_\eta(\alpha)$. Then one would want to compare the accuracy for these two values, $\alpha^*$ and $\alpha^*_\eta$, and the two other measures of performance, top-$k$ and close enough.}
        
        \textcolor{black}{ However, fully identifying the minimizants, $\alpha^*$ and $\alpha^*_\eta$, even in toy settings is a very a difficult task because of non-convexity and dependence on the data-set.} That is why we will consider the following three subsets for some cut-off $\mu>0$ to be determined in terms of $\eta$
          \[
          \begin{split}
&            \mbox{Well classified}=W(\alpha)=\{s,\;\delta \hat X(\alpha,s)>\mu\},\\
&            \mbox{Poorly classified}=P(\alpha)=\{s,\;\delta \hat X(\alpha,s)<-\mu\},\\
&            \mbox{Marginally classified}=M(\alpha)=\{s,\;|\delta \hat X(\alpha,s)|\leq\mu\}.\\
                        \end{split}
          \]
          \textcolor{black}{In our future analysis we plan to address the following issues:
          \begin{itemize}
          \item The minimizant for the 2-scale loss, $\alpha^*_\eta$, will have a low mass (measure) of miss-classified objects whenever the minimizant for the single scale loss,$\alpha^*$, has a low mass of miss-classified objects. In more precise terms, we conjecture the mass (measure) of poorly classified objects $|P(\alpha^*_\eta)|$ for the minimizant $\alpha^*_\eta$ is at worst comparable to the mass $C\,|P(\alpha^*)|$ for some constant $C>0$, depending on the data set and network structure. This ensures that if there are few poorly classified objects for $\alpha^*$ then there are also few poorly classified objects for $\alpha^*_\eta$.
          \item We also conjecture that the total mass (measure) $|P(\alpha^*_\eta)|$ cannot be larger than $C\,|M(\alpha^*_\eta)|$. This is of course not the case for $\alpha^*$ where $|P(\alpha^*)|$ would be much larger than $|M(\alpha^*)|$. In other words, the 2-scale loss function should prevent too many objects from being poorly classified.
            \item We expect that using the 2-scale loss function $L_{2,\eta}(\alpha,R_1,R_2)$, will improve the stability in the sense of \cite{BerJabSaf20}.  This is because we would have a necessarily small loss function when the mass (measure) of well classified objects $|W(\alpha)|$ is close to one; whereas this may not be the case for the classical loss function $L(\alpha)$ (see \cite{BerJabSaf20}).\\
            \end{itemize}
        }

	\section{Numerical Results}
	\label{Sec:Num}
	\subsection{Measures of performance}
	Before comparing the  two-scale loss network with the traditional single scale loss network, we first introduce the measures of performances that we will use to compare the networks.  The first measure of performance is {\it overall accuracy}, referred to as accuracy for convenience, as defined in \eqref{acc_defn}.  Accuracy measures the proportion of objects which are correctly classified, and it can be seen as a good approximation of the likelihood of a new object not used for training being classified correctly by the DNN.
	
	However, accuracy certainly does not give the entire picture about the performance of a classification DNN.  As we mentioned in the introduction, accuracy cannot differentiate between an object that is very poorly classified and objects that are almost well classified. 
	The {\it top-$k$} measure of performance is a commonly used alternative to accuracy which can capture this distinction.  In the top-$k$ measure of performance an object is counted if the correct class has one of the $k$ largest probabilities.  This yields the following definition of top-$k$ performance on the test set $T$,
	\[
	acc_{top-k}=\dfrac{\#\{s\in T: p_{i(s)} \text{ is one of the $k$ largest probabilities}\}}{\#T}, 
	\]
	and a similar definition can be created for the training set by replacing $T$ with $TR$.  For data-sets with a large number of classes or very similar classes, such as oak tree and maple tree in CIFAR100, this is very useful.  For example for CIFAR100 with 100 classes, an object considered as accurate by $acc_{top-2}$ but not correctly classified still has a better predicted probability for the right class larger than 98\% of other predicted probabilities.  
	
	However, a drawback of top-$k$ performance is that an object $s$ can be included as performing well when the actual probability of belonging to the correct class is very low.  For example if the correct class for an object $s$ is the second class, then both the top-$k$ performance and the overall accuracy cannot distinguish between the two cases $(0.9,0.1,0,\dots,0)$ and $(0.51,0.49,0,\dots,0)$.  However it is quite clear that the second case is preferable and the first case is actually even worse than simply assigning every class the same probability. 
	%
	
	
	For this reason we introduce a family of so-called {\it close-enough} measures of performance which counts the proportion of objects within a fixed probability, $at$, of being classified correctly,
	\[
	acc_{at}=\dfrac{\#\{s\in T: p_{i(s)}\geq\max_{1\leq i\leq K}p_i-at\}}{\#T}.
	\]
	The measure of performance $acc_{at}$ counts objects which are close to being correctly classified in terms of probability.  \textcolor{black}{By varying $0\leq at \leq 1$ we can observe the quantity of objects the network has slightly, moderately, or egregiously miss-classified.  This measure of performance allows to observe the distribution on poorly classified objects, which accuracy is unable to do.}  One flaw is that an untrained network can achieve very high performance with this measure as all probabilities are nearly equal. Therefore it is necessary to look at the other measures of performance in conjunction with close-enough performance. We look at all three measures of of performance on CIFAR100 to give a more complete picture of the performance of the DNN.
	
	\subsection{CIFAR100 Numerical comparison: Two-scale loss function with varying scales }
	\label{Sec:Num:2_scale}
	%
	To check the potential improvement in performance numerically, we  consider the first training scheme \eqref{2_scale_loss} from Sec. \ref{Sec:2_loss_var} on a variety of data-sets including MNIST, CIFAR10 and CIFAR100.  For all data-sets a stochastic gradient descent algorithm is implemented with a batch size of 128 and learning rate of $\tau=0.1$.  Various choices of $0<\eta$ were tested.  In addition $R_1$ and $R_2$ were chosen initially as $R_1=R$ and $R_2=10R_1$, where $R$ is the size, \eqref{R_def},  of the initial choice of parameters $\alpha^0$. To train the network with a traditional one scale loss function we simply replace $R_2$ by $R_1$ in \eqref{2_scale_loss}.  In that case we have a single scale traditional loss function, which we use as a baseline to compare the performance of the two-scale loss network with. The python package PyTorch \cite{pytroch} was used to build and train the neural networks as well as the two-scale loss functions.  
	
	{The two-scale loss function was tested on the CIFAR100 data-set, a set of color images with 100 classes (dogs, cars, deer, oak tree, maple tree, etc.) for two types of DNNs, a simple network and a network with state of the art performance, wide residual network (WRN) with adaptive sharpness-aware minimization (ASAM) from \cite{KwoKimPar21}.   The structure of the simple network used is LeNet-5 adjusted for CIFAR100 \cite{LecJacBot}.}  For the simple network there are two sets of convolutional layers, each followed by a max pooling layer.  Afterwards two fully connected layers are used. Lastly, the soft-max function is used to obtain the probabilities as in \eqref{softmax}.  Training was performed for 128 epochs with a batch size of 128,  (one epoch contains $\{\#TR\}/$(batch size) number of steps). For CIFAR100, the networks were tested for 10 different seeds, leading to different initial random $\alpha^0$. The measures of performance were then averaged over all seeds.  For all data-sets, accuracy for the batch used is recorded at each iteration and used to approximate the accuracy of the training set.  
	The accuracy on the test set is measured exactly at the end of training and intermittently during training, with frequency depending on the length of training.  The other measures of performance were measured at the end of training over the entire test set or the entire training set.  Achieving a high level of performance on CIFAR100 is typically more challenging for DNNs since there are relatively few objects (500) per class to train from and a large number of classes \cite{hacdap19}.  In addition several of the classes are similar, e.g. oak tree and maple tree.

	Note that when we compute the probabilities for the test set, we must choose a scale for the two-scale loss network as the network should not use the test set's objects' correct classes in order to classify the objects or assign probabilities.  The correct class should only be used to determine performance after the probabilities are determined.  We choose the scale as $R_1$, as picking $R_2$ would yield very large and small probabilities. 
	
\begin{figure}[!t]
	\centering
	
	\includegraphics[scale=.5]{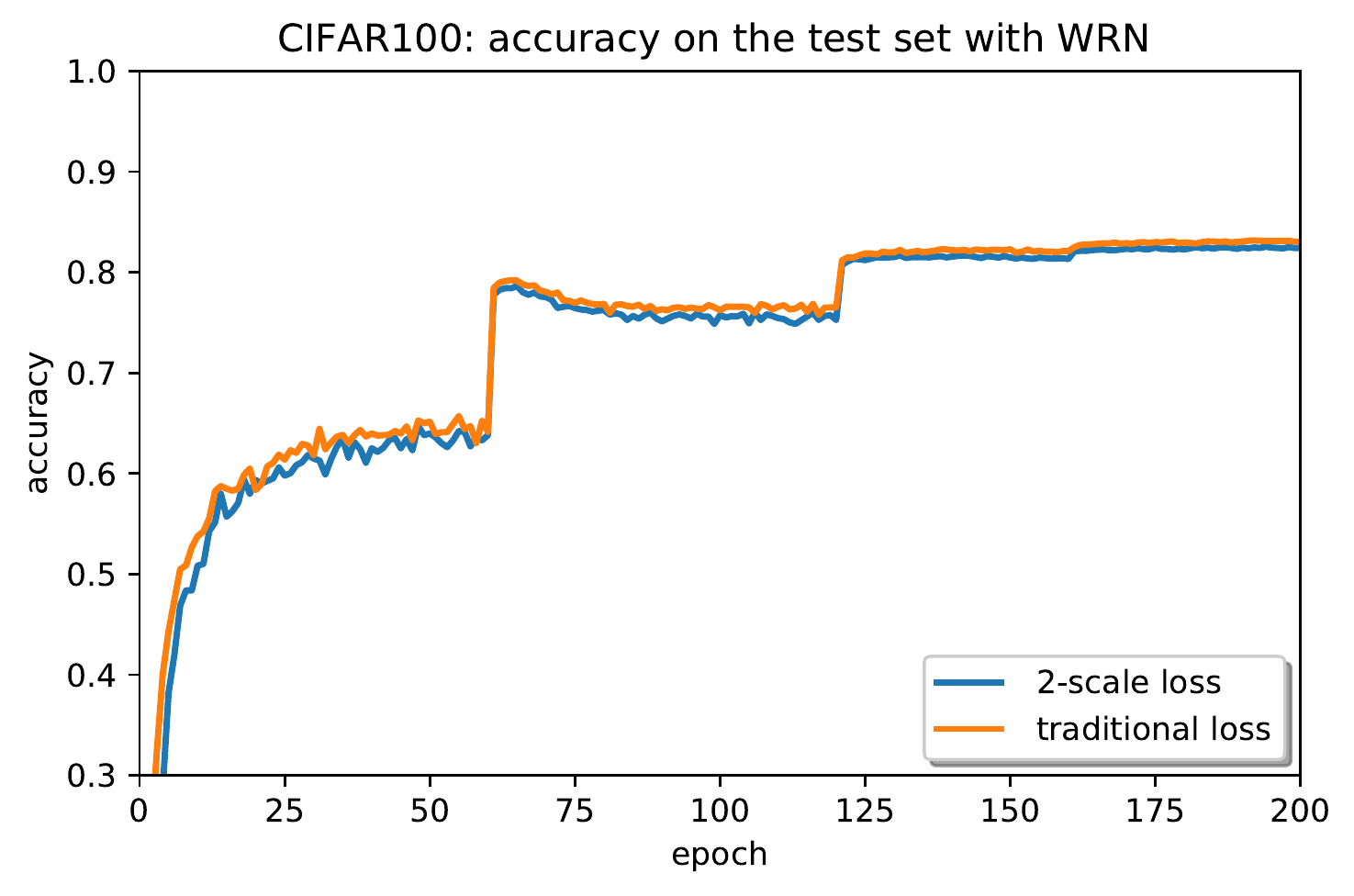}
	\includegraphics[scale=.5]{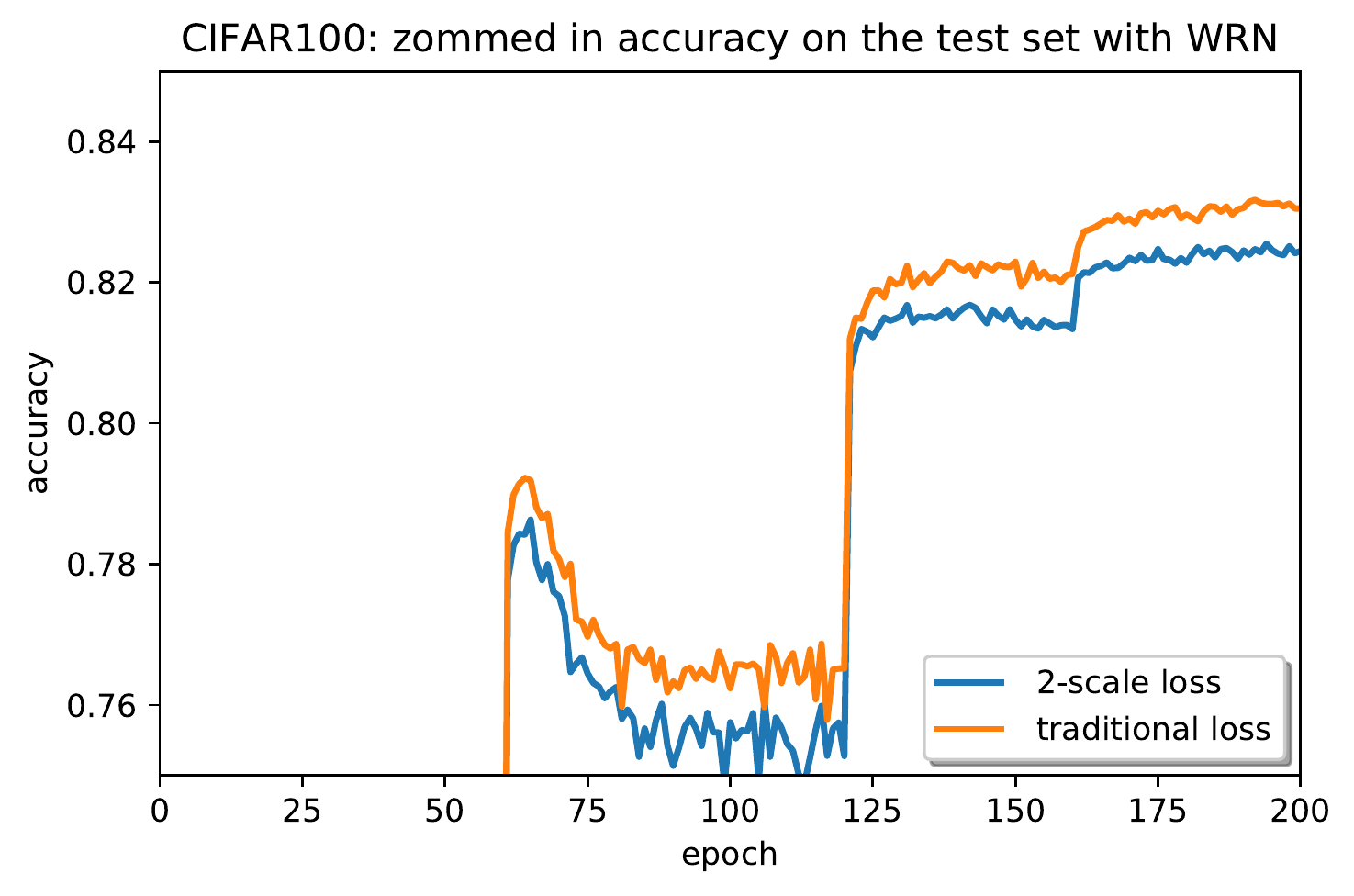}
	\caption{{The accuracy throughout training on the test set for two-scale loss WRN \eqref{2_scale_loss} was compared with the traditional single scale loss WRN for various values of $\eta$ on the CIFAR100 data-set.  The performance on the test set was recorded approximately every 390 iterations ( 1 epoch) and the graphs represent the average accuracy of 5 random initial seeds.  The right graph is a zoomed in version of the graph on the left. \label{Fig:2_scale_CIFAR100_sam}}}
\end{figure}
\begin{figure}[!t]
		\centering
		
		\includegraphics[scale=.5]{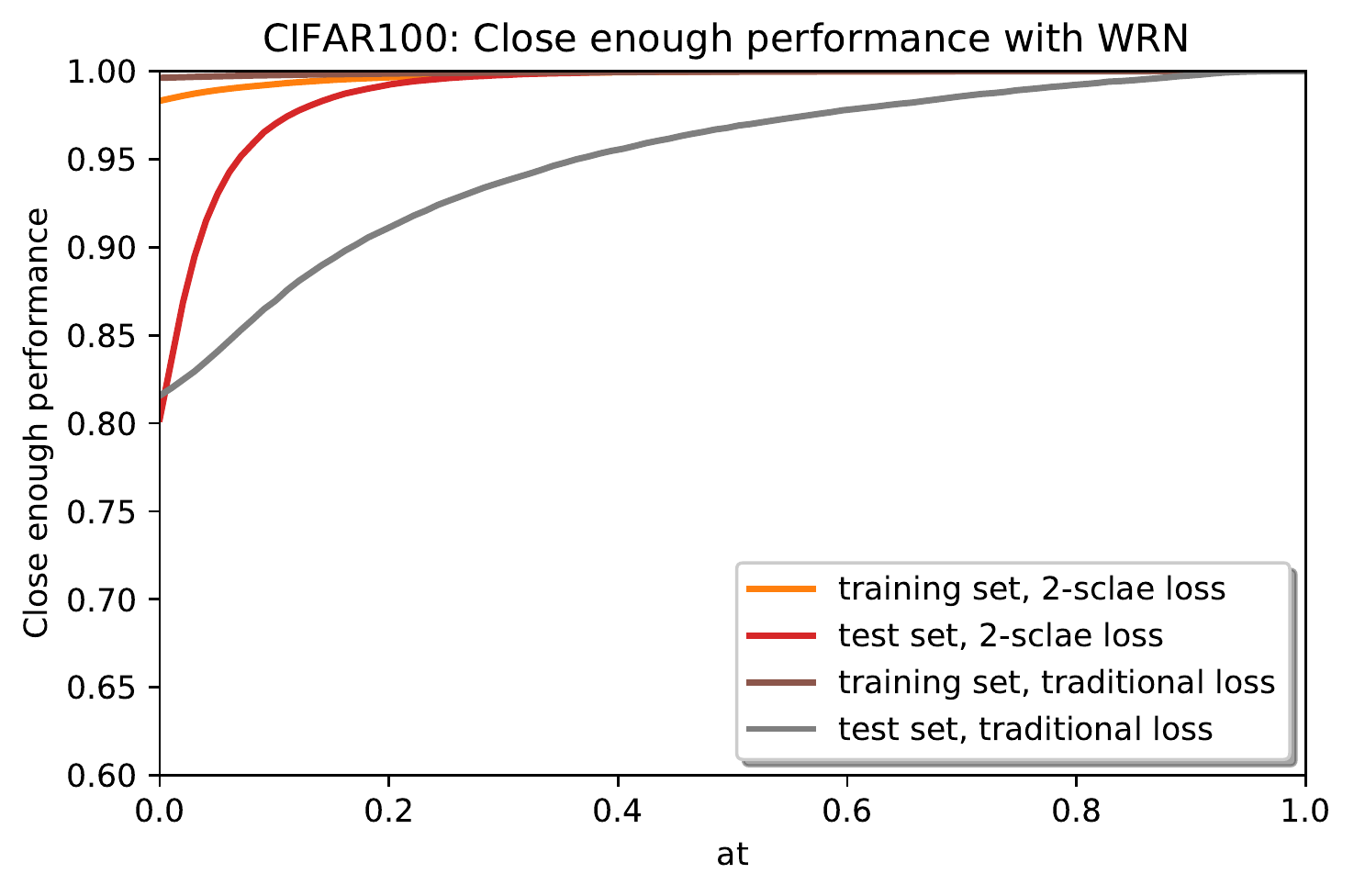}
		\caption{{The accuracy throughout training on the test set for two-scale loss WRN \eqref{2_scale_loss} was compared with the traditional single scale loss WRN for various values of $\eta$ on the CIFAR100 data-set.  The performance on the test and training set was recorded at the end of training( 200 epochs) and the graphs represent the average performance of 5 random initial seeds. \label{Fig:close_enough_CIFAR100_sam}}}
	\end{figure}

{	First we compare the performance on the WRN with ASAM.  We observe fairly high accuracy however the two-scale loss performs slightly worse at the end of training with roughly 0.5\% lower accuracy on the test set, see Fig. \ref{Fig:2_scale_CIFAR100_sam}.  However, as mentioned previously, accuracy does not discern how the network performs on miss-classified objects which make up a significant portion of objects in the test set at the end of training. Therefore we next look at the close-enough performance where we observe a large increase in close-enough performance on the test set for the two-scale loss network for values of $0.01<at<0.8$ with a maximum increase of 10\% when $at=0.1$ , see Fig. \ref{Fig:close_enough_CIFAR100_sam}.  This large increase of performance shows that the two-scale loss is ensuring that the network performs much better on the miss-classified objects.  For the traditional loss, there are many objects which are far from begin correctly classified while the two-scale loss network has almost all objects counted in the close-enough performance by $at\approx0.1   $.}
	
		\begin{figure}[!t]
		\centering
		
		\includegraphics[scale=.5]{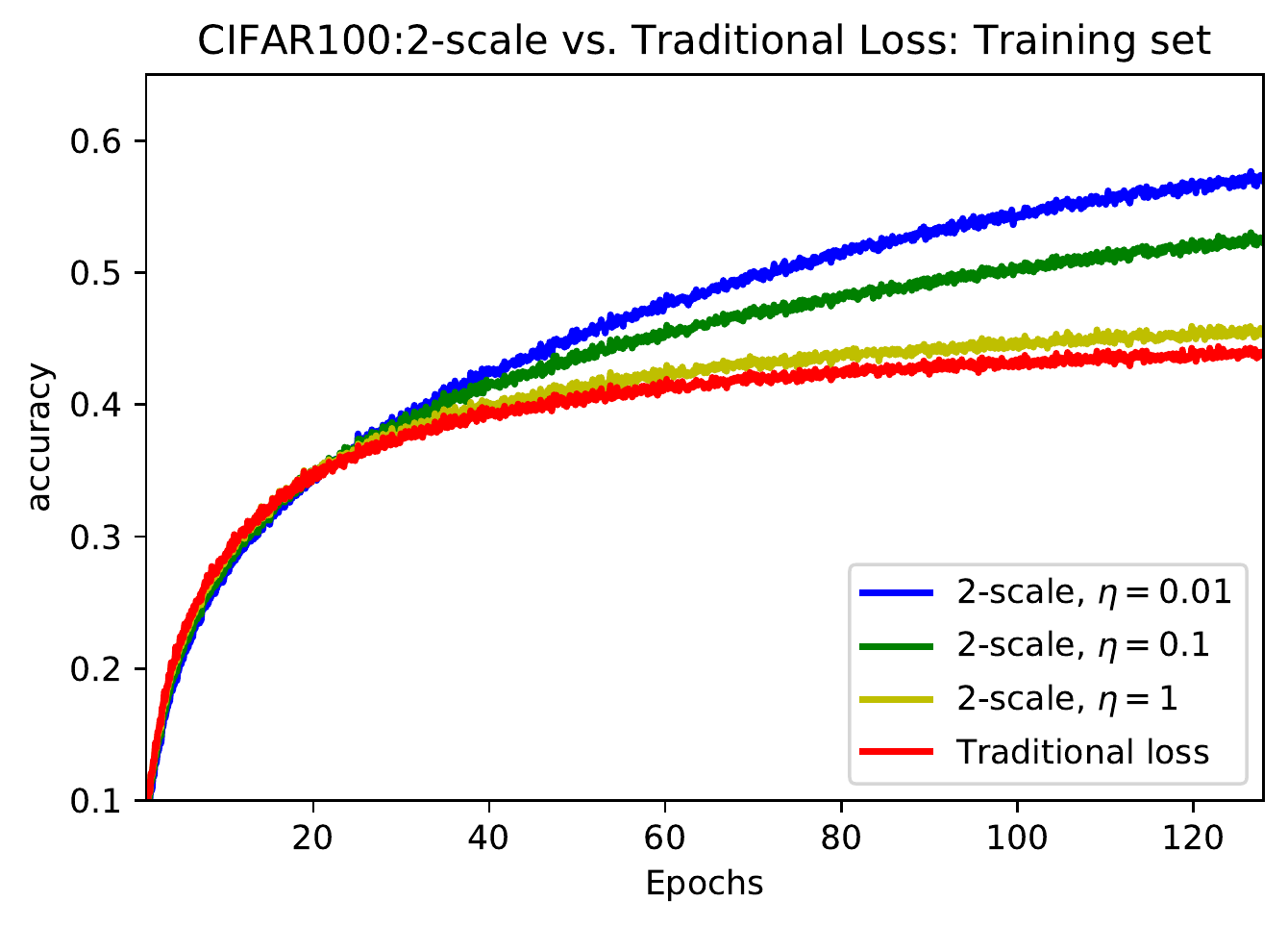}
		\includegraphics[scale=.5]{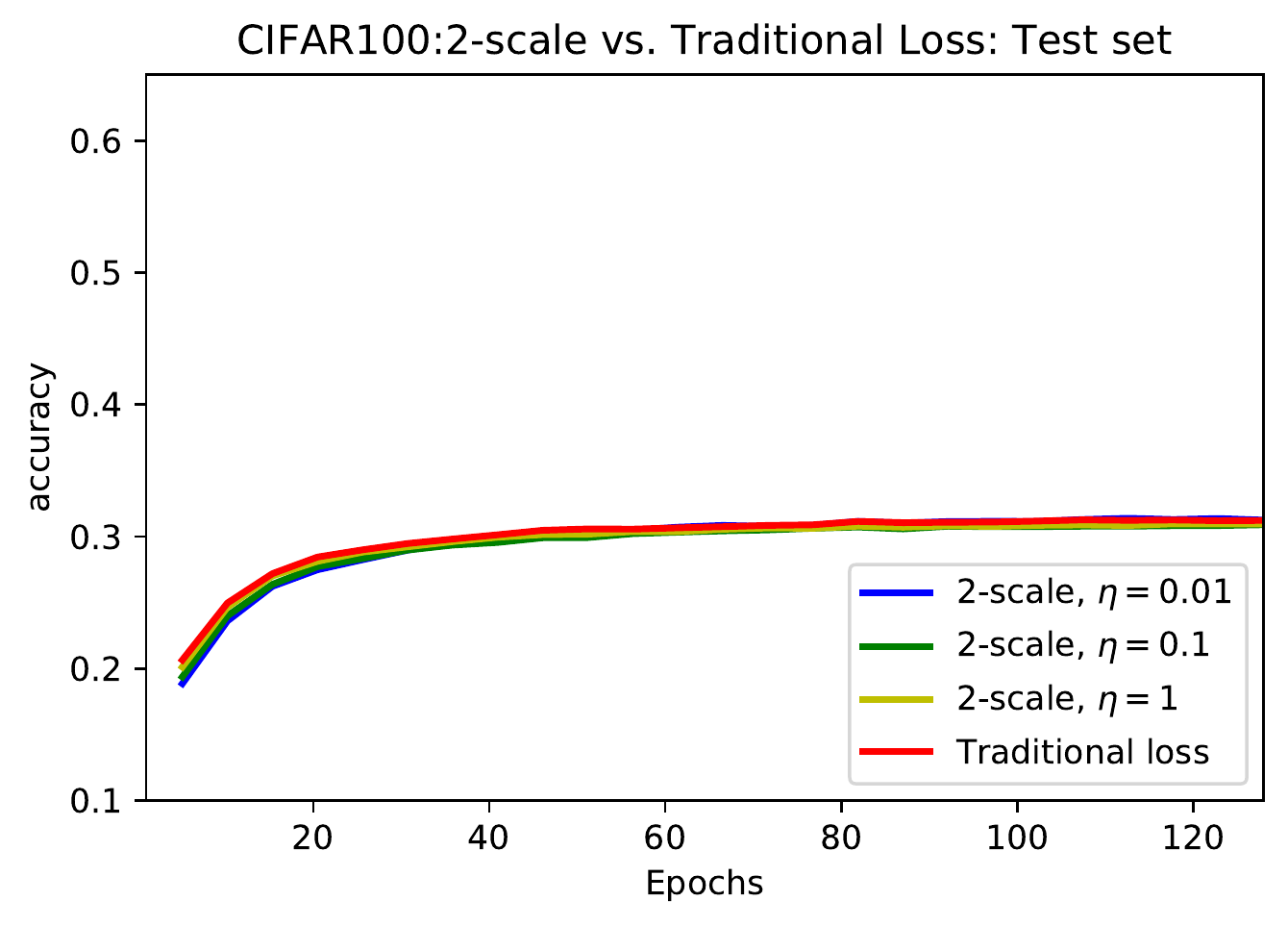}
		\caption{The accuracy throughout training on the training set (left) and the test set (right) for two-scale loss network \eqref{2_scale_loss} was compared with the traditional single scale loss network for various values of $\eta$ on the CIFAR100 data-set.  Accuracy on the test set was recorded once every 2000 iterations ($\approx$5 epochs) and the graphs represent the average accuracy of 10 random initial seeds. \label{Fig:2_scale_CIFAR100_lt}}
	\end{figure}

	\begin{figure}[!t]
		\centering
		\includegraphics[scale=.5]{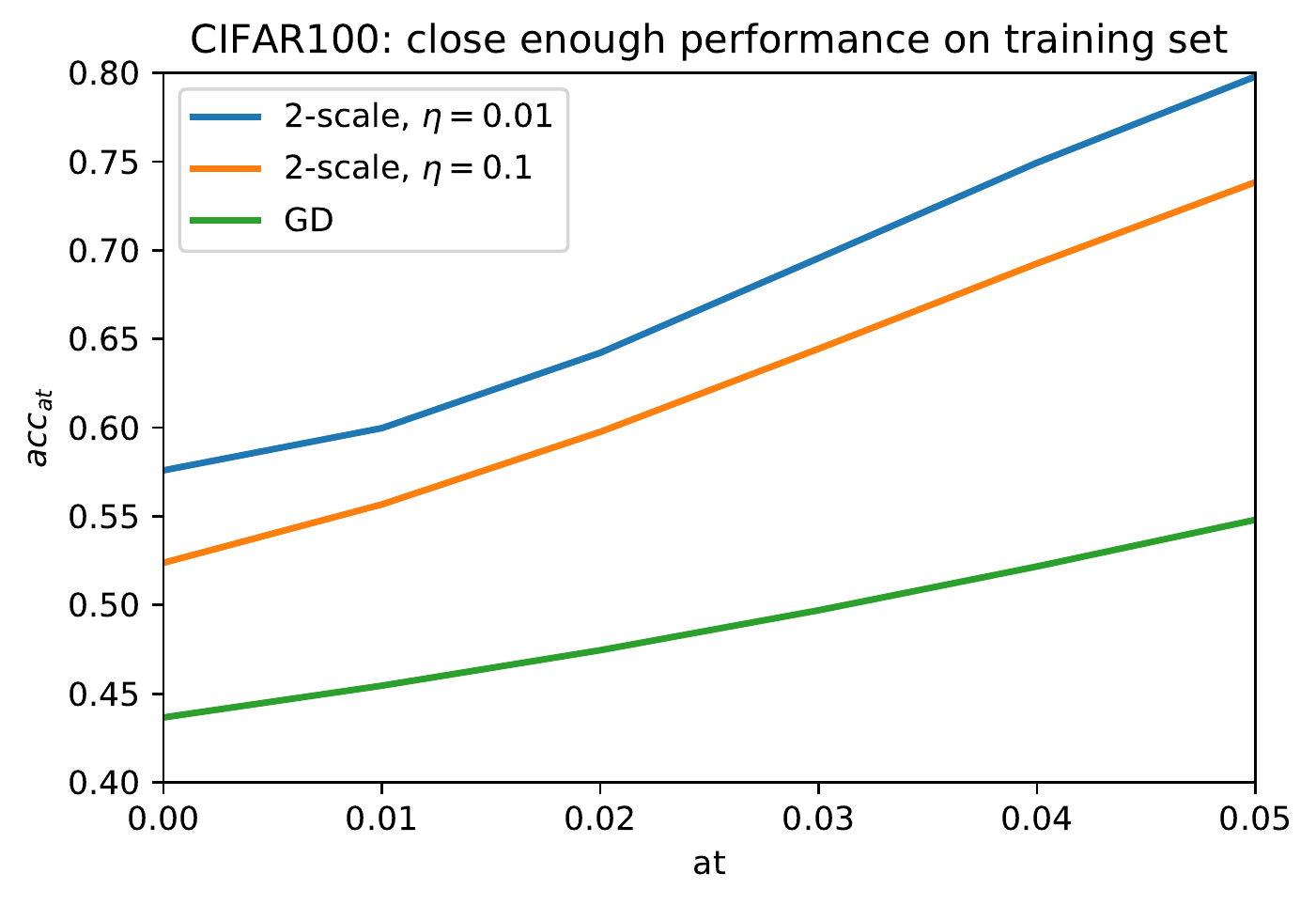}
		\includegraphics[scale=.5]{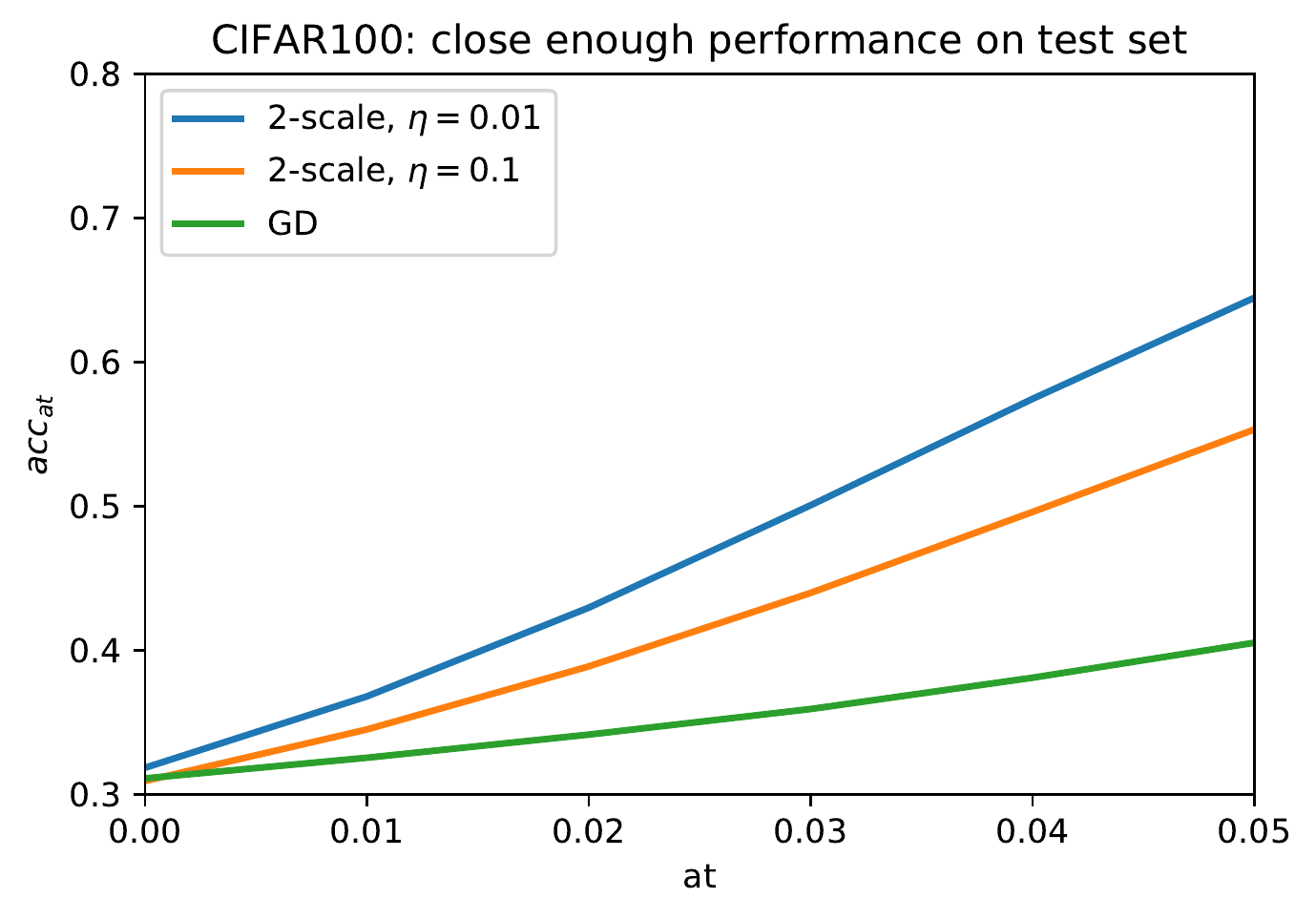}
		\caption{The close-enough performance on the training set (left) and the test set (right) for the two-scale loss network \eqref{2_scale_loss} was compared with the traditional single scale loss network for various values of $\eta$ on the CIFAR100 data-set.  The performance is measured after 128 epochs of training.  The graphs represent the average performance of 10 random initial seeds.   \label{Fig:2_scale_CIFAR100_at}}
	\end{figure}
	
	\begin{figure}[!t]
		\centering
		\includegraphics[scale=.5]{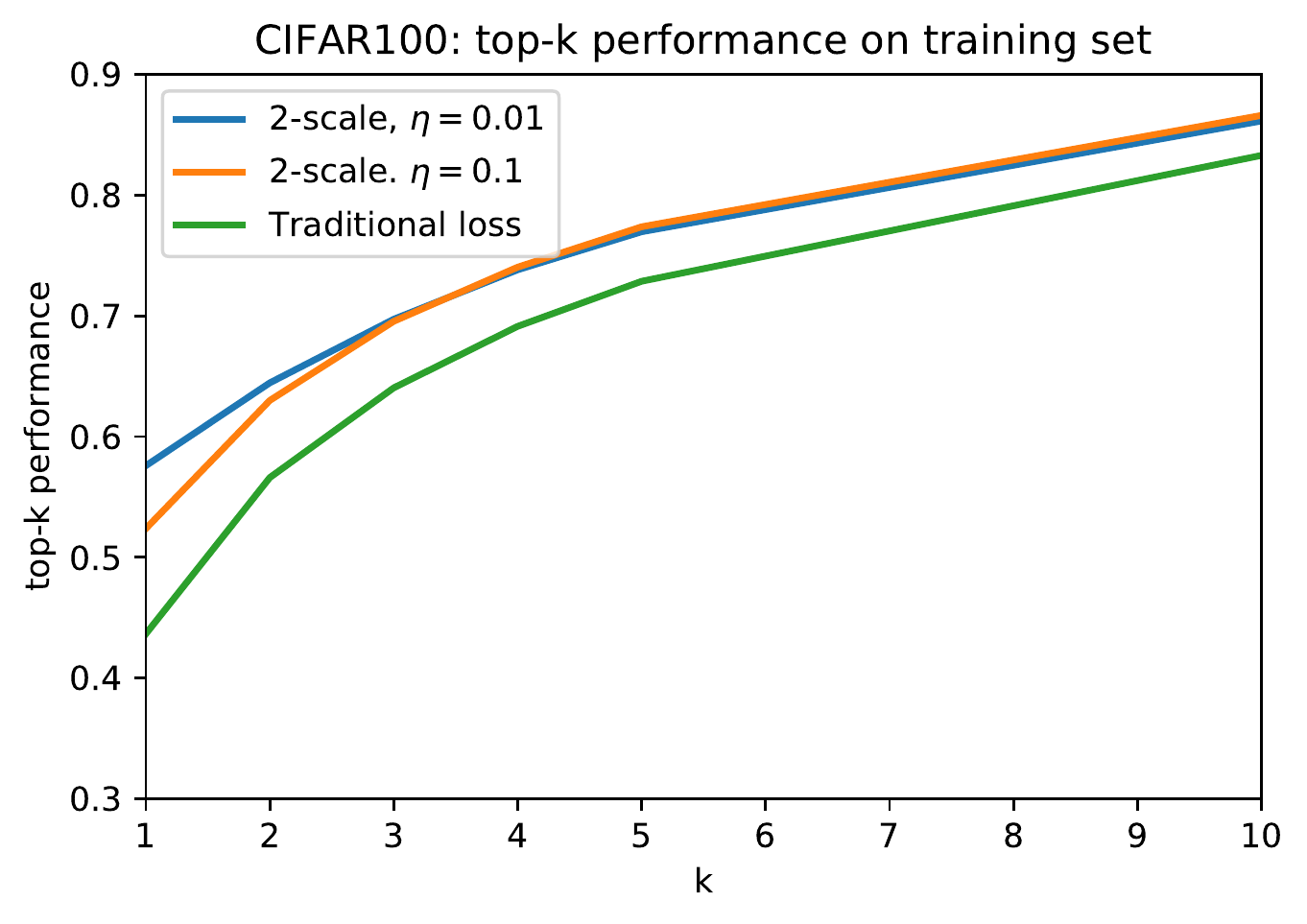}
		\includegraphics[scale=.5]{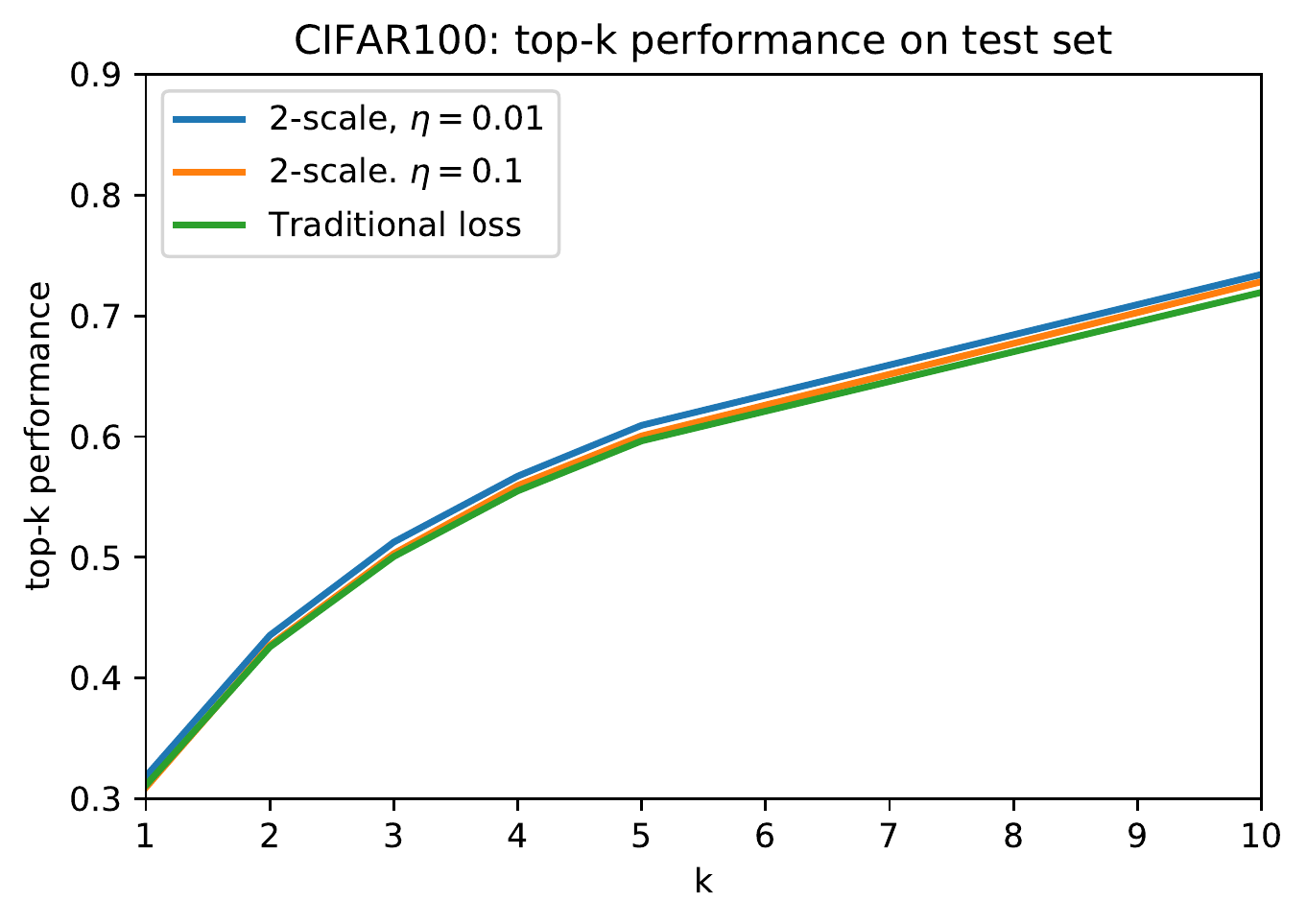}
		\caption{The top-$k$ performance on the training set(left) and the test set (right) for the two-scale loss network \eqref{2_scale_loss} was compared with the traditional single scale loss network for various values of $\eta$ on the CIFAR100 data-set after 128 epochs of training.  The graphs represent the average performance of 10 random initial seeds. \label{Fig:2_scale_CIFAR100_topk}}
	\end{figure}

	Next we compare the performance on the simple network. We begin by observing a large increase in performance for all three measures of performance on the training set by using the two-scale loss network with $\eta=0.01$, see Figs. \ref{Fig:2_scale_CIFAR100_lt}-\ref{Fig:2_scale_CIFAR100_topk}. The two-scale networks have a lower accuracy than the traditional loss up through 20 epochs, and then reverse the trend and attain a much larger accuracy, see Fig \ref{Fig:2_scale_CIFAR100_lt}.  After 128 epochs, with $\eta=.01$, the two-scale network has an average increase of accuracy of 13.4\%. We also observe as that $\eta$ increases the accuracy of the two-scale network approaches the traditional loss function, as expected since for larger $\eta$ the second scale will be used less.

	At first look, the performance on the testing set appears indistinguishable for the two-scale loss and the classical cross-entropy  with the two-scale network having an average increase of accuracy of 0.1\% near the end of the training. While this is correct when looking only at the overall accuracy, it masks an actual large increase in performance on not well classified objects. This point is made obvious when looking at the close-enough performance, $acc_{at}$, on the test set for $0.01<at<0.05$, see Fig \ref{Fig:2_scale_CIFAR100_at}: The two-scale loss network achieves a substantial improvement over the traditional loss network for all $at\geq0.01$ and the improvement increases as $at$ grows, see Fig \ref{Fig:2_scale_CIFAR100_at} (right).   This indicates that the two-scale loss function has still drastically improved  the number of objects that are close to being correctly classified.
	
	%

	Rather naturally, an increase of top-$k$ performance of 1.5\% at $k=10$ on the test set is also observed at the end of training for the two-scale loss network with $\eta=0.01$ and the increase grows slightly larger at $k$ increases, see Fig \ref{Fig:2_scale_CIFAR100_topk} (right).  We want to emphasize that such a 1.5\% increase (3\% in relative performance) comes from a very small increase in computational complexity of a net increase of one parameter in the loss function for a network with over 50,000 parameters.  
	Obtaining any noticeable improvement while only increasing the total parameters by one is a significant gain where improvements of even 0.5\% are often considered significant.  
	When comparing top-$k$ performance the improvement on the training and test set are quite similar for large $k$.
	But  for the training set, we observe that the improvement in top-$k$ performance for two-scale over traditional loss decreases as $k$ increases.
	For small $k$, we again found a larger improvement in performance for the two-scale loss on the training set vs. the test set. 
	

	\begin{figure}[!t]
		\centering
		
		\includegraphics[scale=.5]{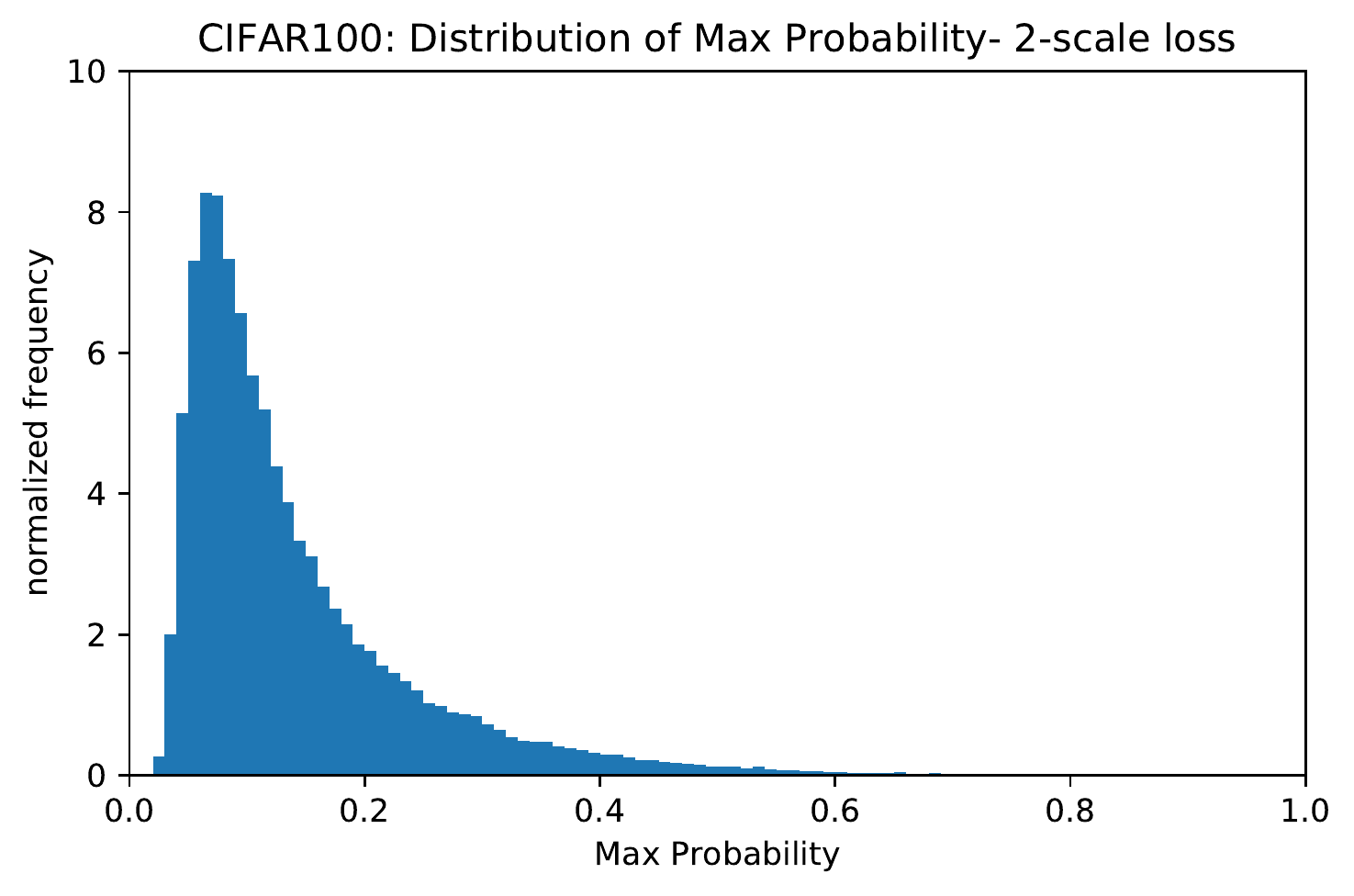}
		\includegraphics[scale=.5]{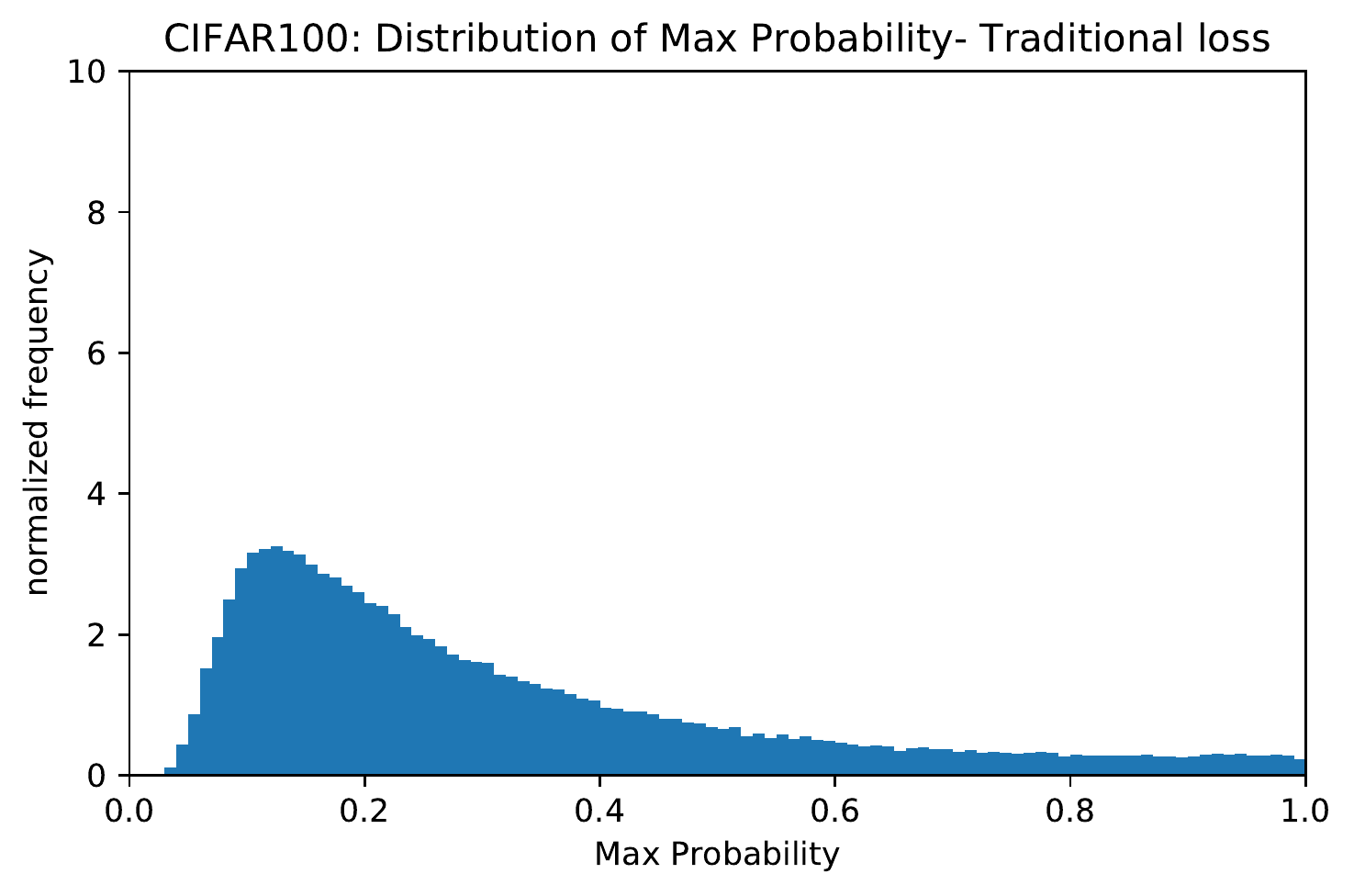}
		\caption{A histogram of the max probabilities of objects on the test set for the two-scale loss network \eqref{2_scale_loss} with $\eta=0.01$ (left) and traditional single scale loss network (right) on the CIFAR100 data-set after 128 epochs and across 10 seeds. \label{Fig:max_prob}}
	\end{figure}	
	
	To understand why see such a large increase in the close-enough performance on the test set in Fig. \ref{Fig:2_scale_CIFAR100_at}, we take a closer look at the distribution of probabilities over the test set.  We first consider the maximal probability assigned for each object on the test set by the DNNs. We immediately emphasize that this max probability may not be the correct probability of course if the object is not well classified.
	
	We observe that the traditional single scale loss network tends to have many more objects with large maximal probabilities in general, see Fig. \ref{Fig:max_prob}.  Some of this effect can just be explained by our changes in the scales: We are using the scale $R_1$ for the two-scale loss and $R_1$ is slightly smaller than the $R$ obtained in the traditional loss.  An increase in $R$, will either increase or maintain the max probability of each object, therefore since $R_1$ is less than $R$ in the traditional single scale loss we expect lower max probabilities for the two-scale loss network.  However, as the accuracy is only approximately 30\%, some of these very large max probabilities are actually from the incorrect classes.  If the maximal probability is assigned to the wrong class, then it is actually beneficial to reduce it. We hence conjecture that the reason the single scale loss function is performing much worse on those misclassified objects is because it assigns them large probabilities on the wrong class.
	
	In addition there are very few objects whose max probability is less than 5\% and none whose max probability is less than 2\% for both networks.  Therefore $acc_{at}$ can indeed be used as a reasonable measure of performance as few objects are automatically included for the values of $at$ tested, $0<at<0.05$, and none are automatically included for $at<0.02$.

	\begin{figure}[!t]
		\centering
		
		\includegraphics[scale=.5]{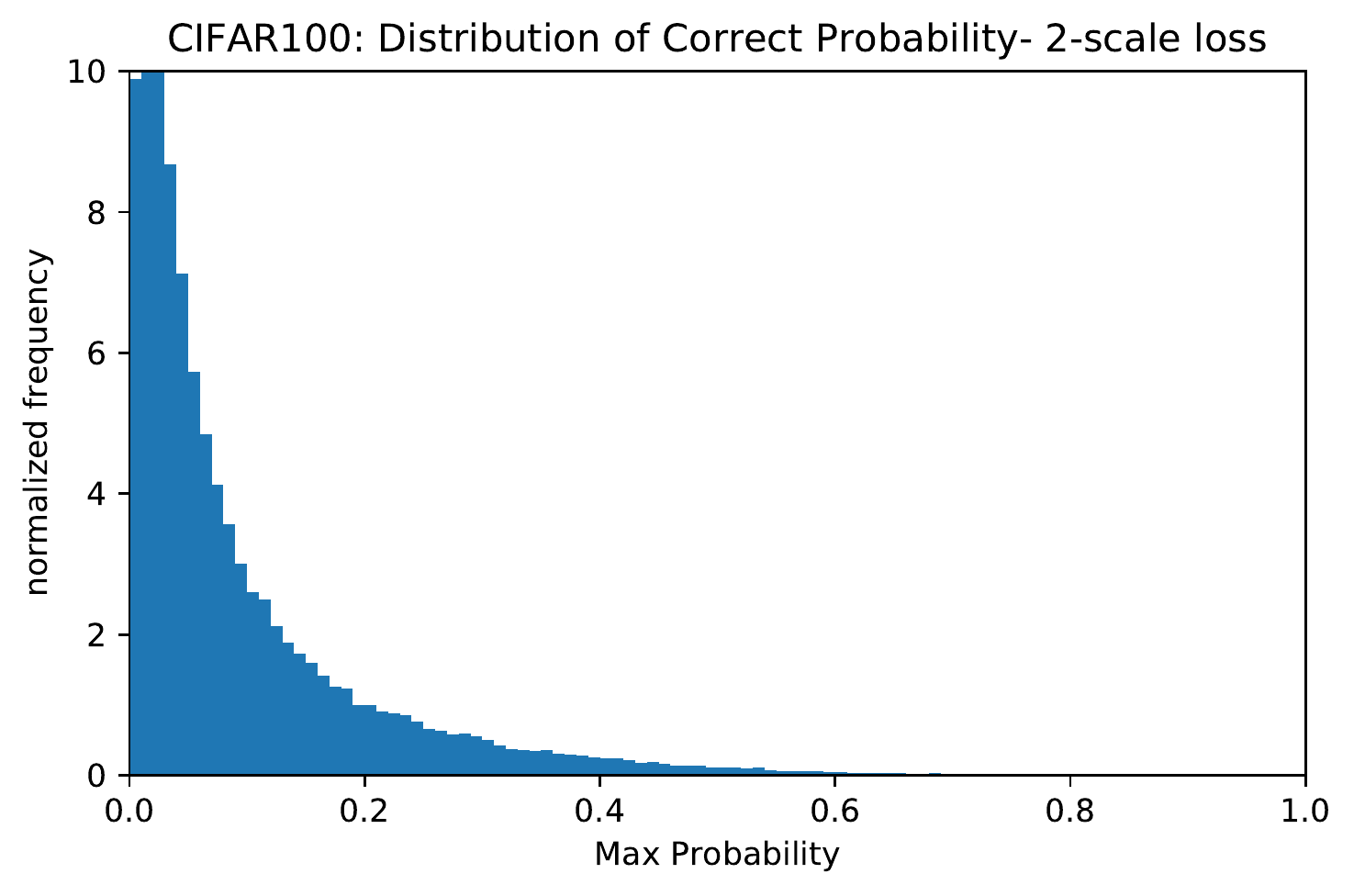}
		\includegraphics[scale=.5]{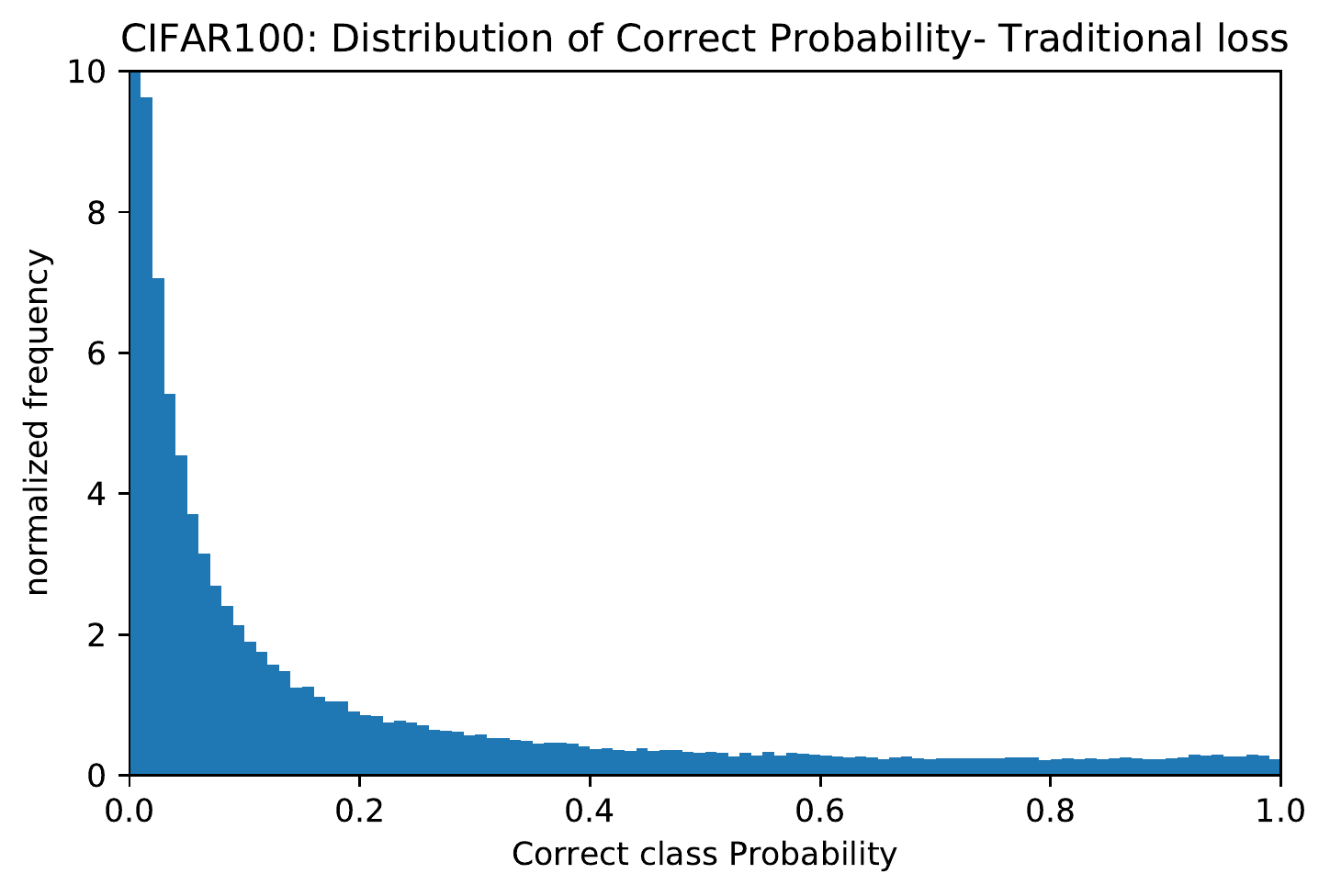}
		\caption{A histogram of the correct probabilities of objects on the test set for the two-scale loss network \eqref{2_scale_loss} with $\eta=0.01$ (left) and traditional single scale loss network (right) on the CIFAR100 data-set after 128 epochs and across 10 seeds. \label{Fig:correct_prob}}
	\end{figure}
	
	We looked at the distribution of probabilities for the correct class of objects, $p_{i(s)}$ on the test set.  We again observe that the single scale loss has more objects with larger probabilities in general, see Fig. \ref{Fig:correct_prob}. This is natural as again $R$ is larger than $R_1$. Moreover this does not tell us how close an object is to being correctly classified or even if the object is correctly classified at all if $p_{i(s)}<0.5$.  For example if an object had $p_{i(s)}=0.2$ and is correctly classified, then we would say this is a better performance than if the object had a slightly larger correct probability, $p_{i(s)}=0.21$, but was miss-classified. 	
	
	For this reason, the histogram of the probability confidence is more relevant to our analysis: the difference between the correct class and the max probability for miss-classified objects,
	\begin{equation}
		\delta p(\alpha,s):=p_{i(s)}(\alpha,s)-\max_{1\leq i\leq K}p_i(\alpha,s).    
	\end{equation}
	\begin{figure}[!t]
		\centering		\includegraphics[scale=.5]{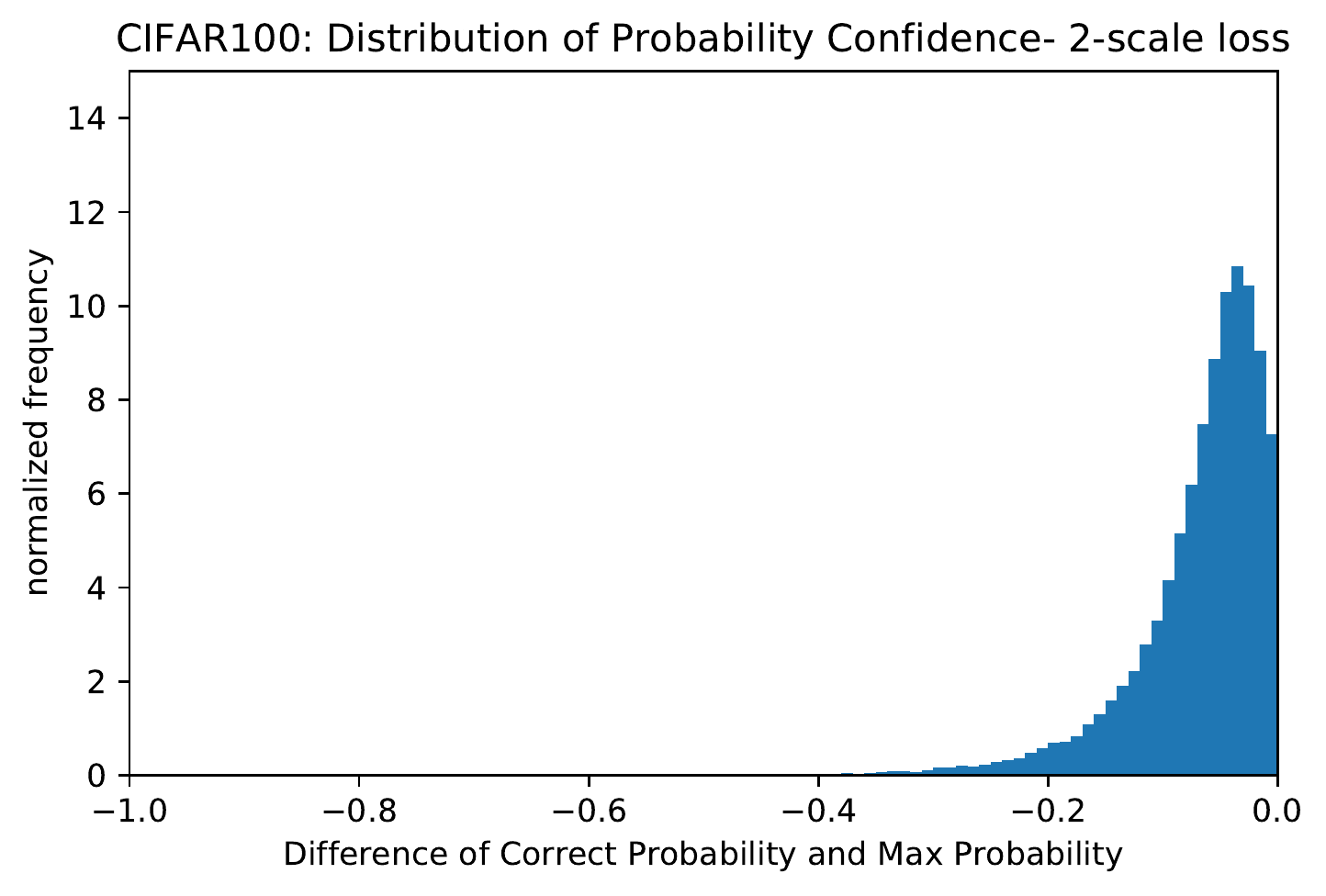}
		\includegraphics[scale=.5]{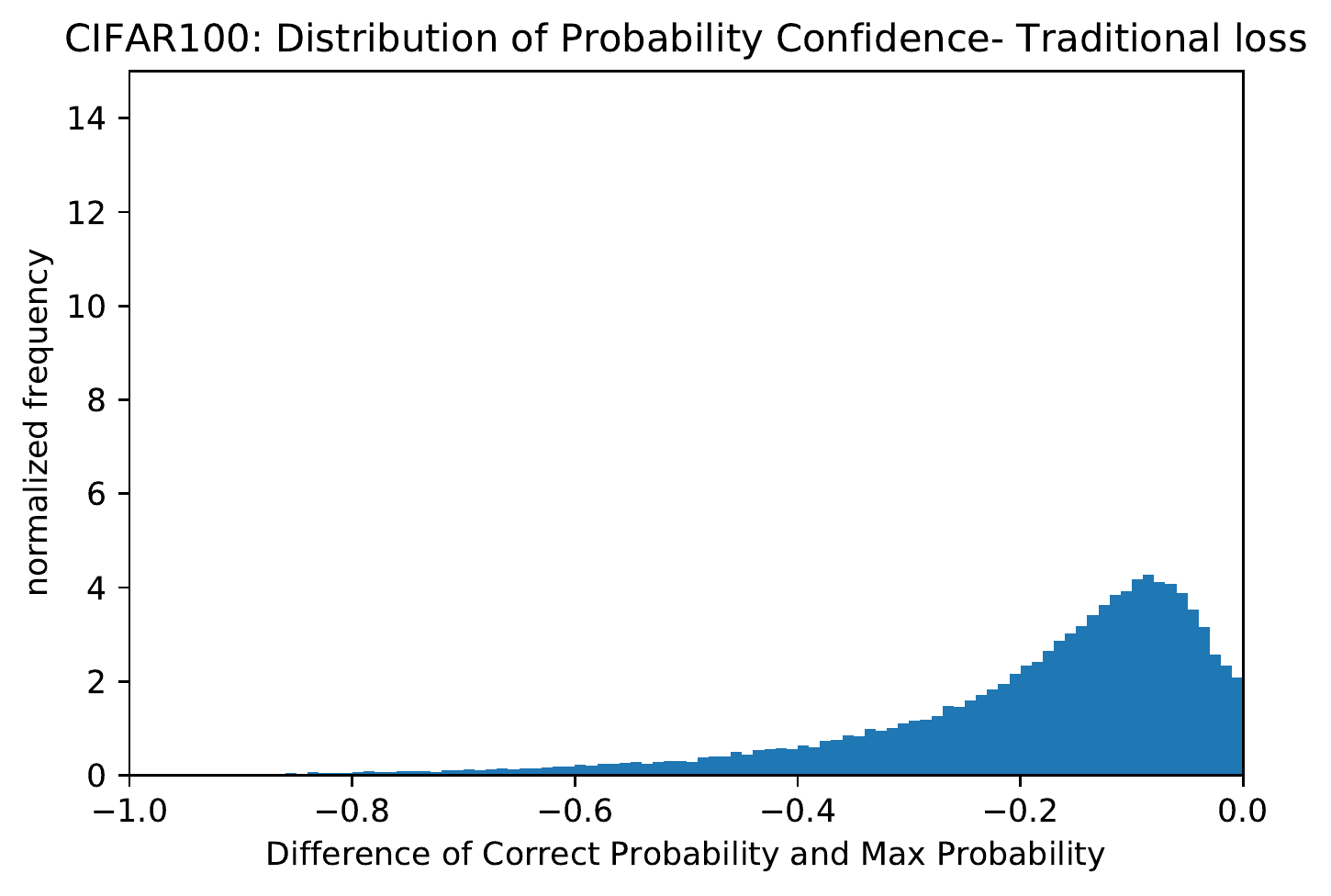}
		\caption{A histogram of the difference of the correct probabilities and max probabilities of miss-classified objects on the test set for the two-scale loss network \eqref{2_scale_loss} with $\eta=0.01$ (left) and traditional single scale loss network (right) on the CIFAR100 data-set after 128 epochs and across 10 seeds. \label{Fig:diff_prob}}
	\end{figure}
	We observe, for the traditional loss network that many objects are miss-classified by a large probability while for two-scale loss network most miss-classified objects are miss-classified by a small probability, see Fig. \ref{Fig:diff_prob}.  This confirms the conjecture that the two-scale loss network is performing much better on objects which are miss-classified by decreasing the difference between the probability on the correct class and the higher probability on some incorrect class(es).

	Finally, we also investigated the relative performance of the two algorithms on CIFAR100 so-called super-classes. In addition to the original 100 classes, CIFAR100 also has 20 super-classes where each super-class is composed of 5 similar classes of the original 100 (e.g. the tree super-class is composed of the maple, oak, pine, willow, and palm classes) \cite{KriHin09}.  Measuring the super-class accuracy as
	\begin{equation}
		acc_{super-class}=\dfrac{\#\{s\in T: \text{ the index for the maximum $p_i$ belongs to the same super-class as $i(s)$}\}}{\#T},
	\end{equation}
	takes into account the similarity of classes, which the other measures of performance we  have used are unable to do.  Of course the super-class accuracy also ignores the probability of belonging the original correct class.

	\begin{figure}[!t]
		\centering
		\includegraphics[scale=.5]{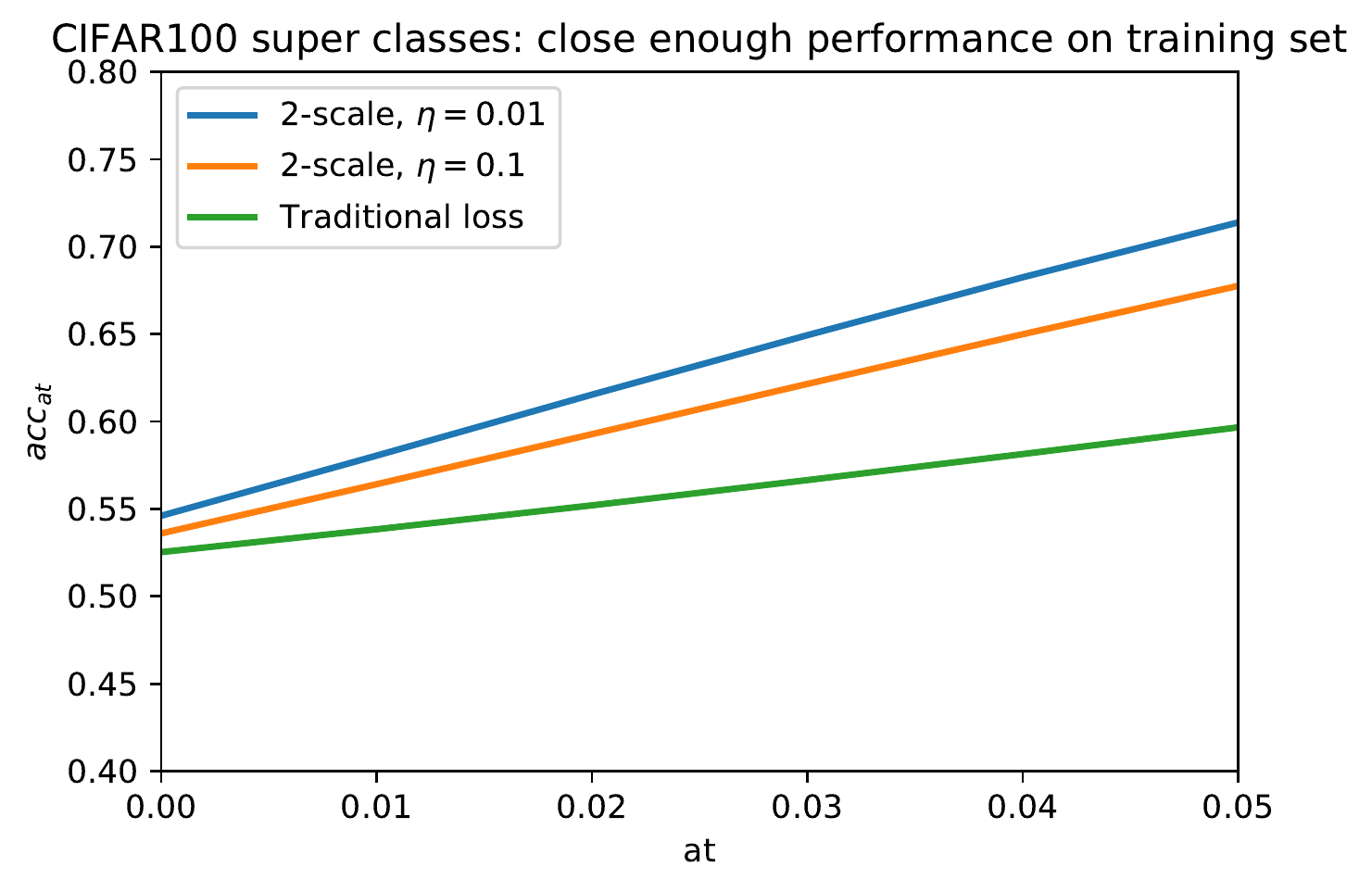}
		\includegraphics[scale=.5]{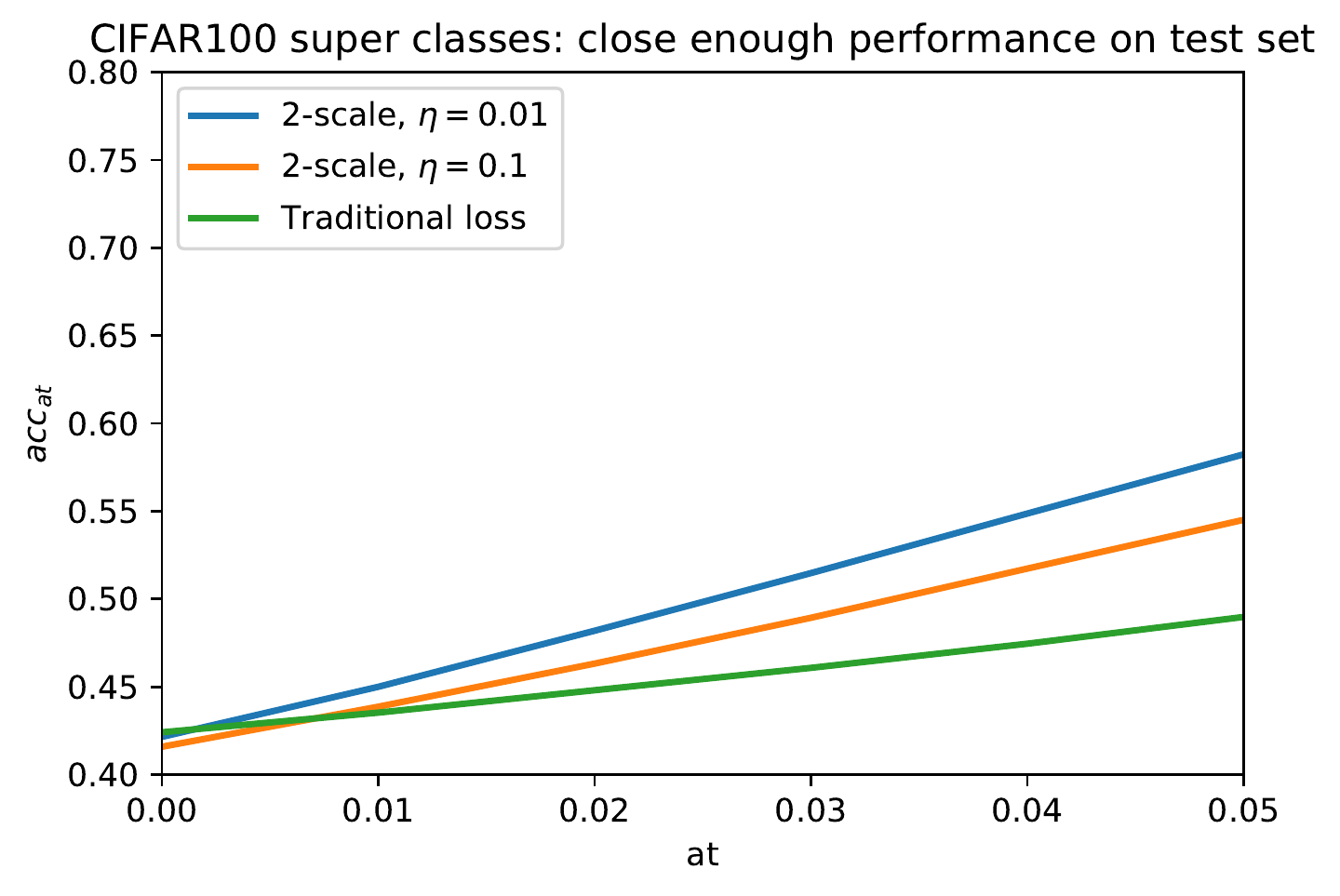}
		\caption{The close-enough performance on the training set (left) and the test set (right) for the two-scale loss network \eqref{2_scale_loss} was compared with the traditional single scale loss network for various values of $\eta$ on the CIFAR100 data-set's super-classes.  The performance is measured after 128 epochs of training.  The graphs represent the average performance of 10 random initial seeds.   \label{Fig:2_scale_CIFAR100_coarse}}
	\end{figure}

	The results for super-classes turn out to be quite similar to the previous discussion. We see a drastic improvement in close-enough performance for the two-scale loss function both on the training and testing sets, see Fig. \ref{Fig:2_scale_CIFAR100_coarse}. In terms of overall accuracy, the two-scale network enjoys a significant gain on the training set, but very similar results on the testing set.	

	The conclusion of this analysis, based on three different measures of performance, is that  the two-scale loss network has many more objects on the test set which are close to being correct with regards to both
	\begin{itemize}
		\item[i] how many incorrect probabilities are higher than the correct one, $p_{i(s)}$,
		\item[ii] how close $p_{i(s)}$ is to the maximum $p_i$.
	\end{itemize}
	This justifies the original motivation of focusing training on not-well classified objects in order to improve their performance.

	\subsection{CIFAR10 and MNIST Numerical comparison: Two-scale loss function with varying scales} 

	The two-scale loss network was next tested on the CIFAR10 data-set \cite{KriHin09}, a set of color images with 10 classes (dogs, cars, deer, etc.).  The structure of the network used is LeNet-5 adjusted for CIFAR10 \cite{LecJacBot}.  There are two sets of convolutional layers, each followed by a max pooling layer.  Afterwards two fully connected layers are used. Lastly, the soft-max function was used to obtain the probabilities as in \eqref{softmax}.  
	
	
	\begin{figure}[!t]
		\centering
		\includegraphics[scale=.5]{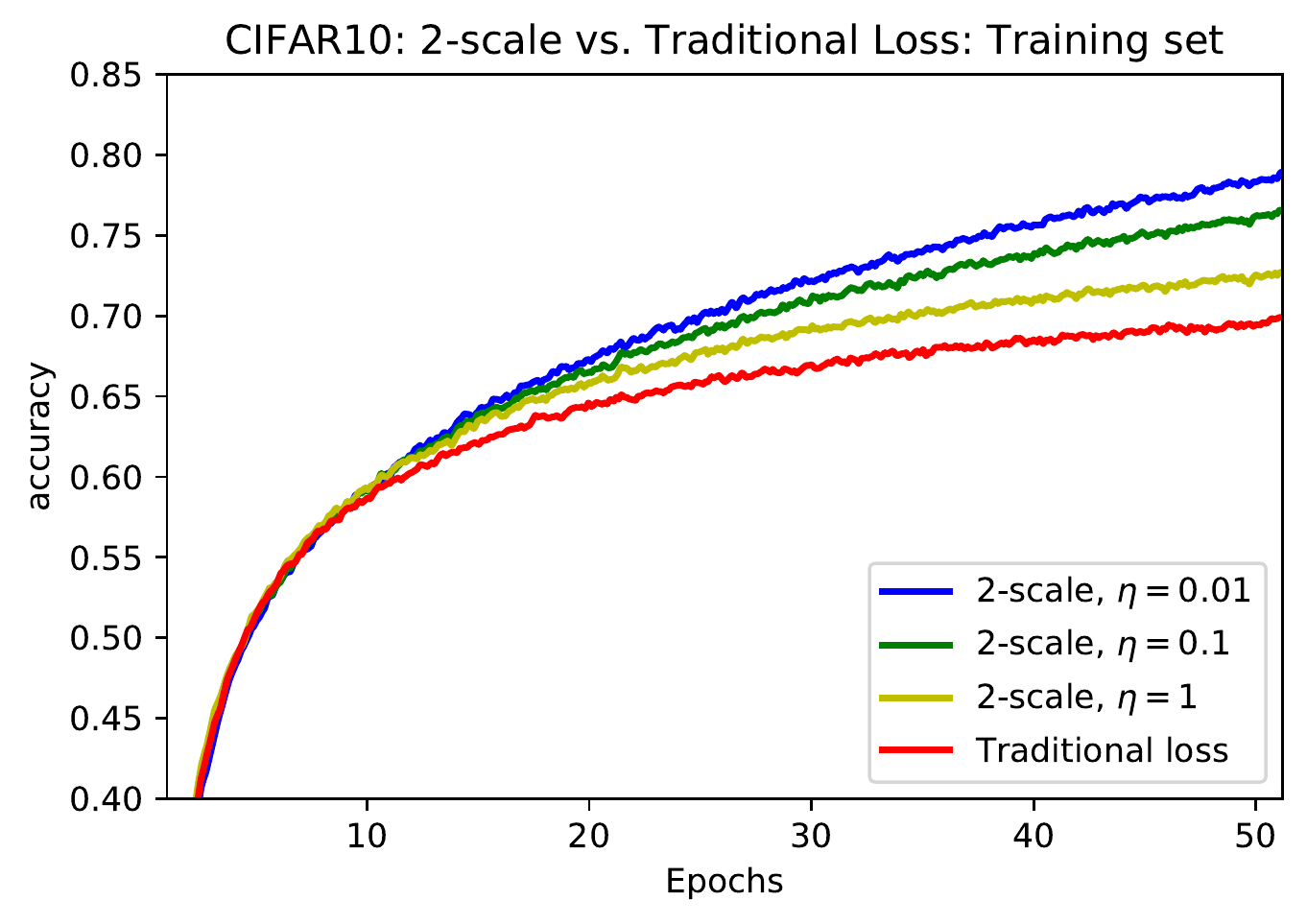}
		\includegraphics[scale=.5]{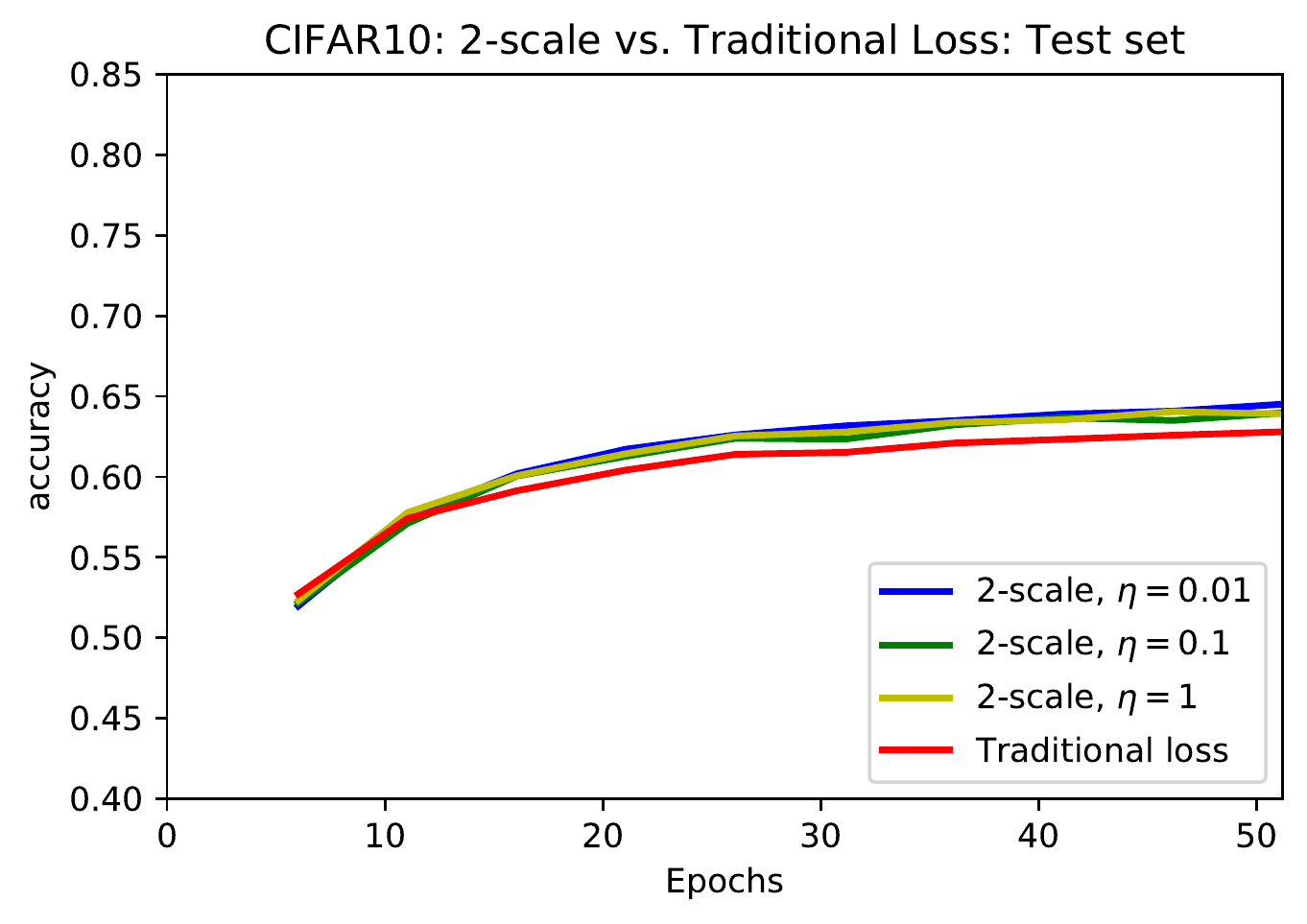}
		\includegraphics[scale=.5]{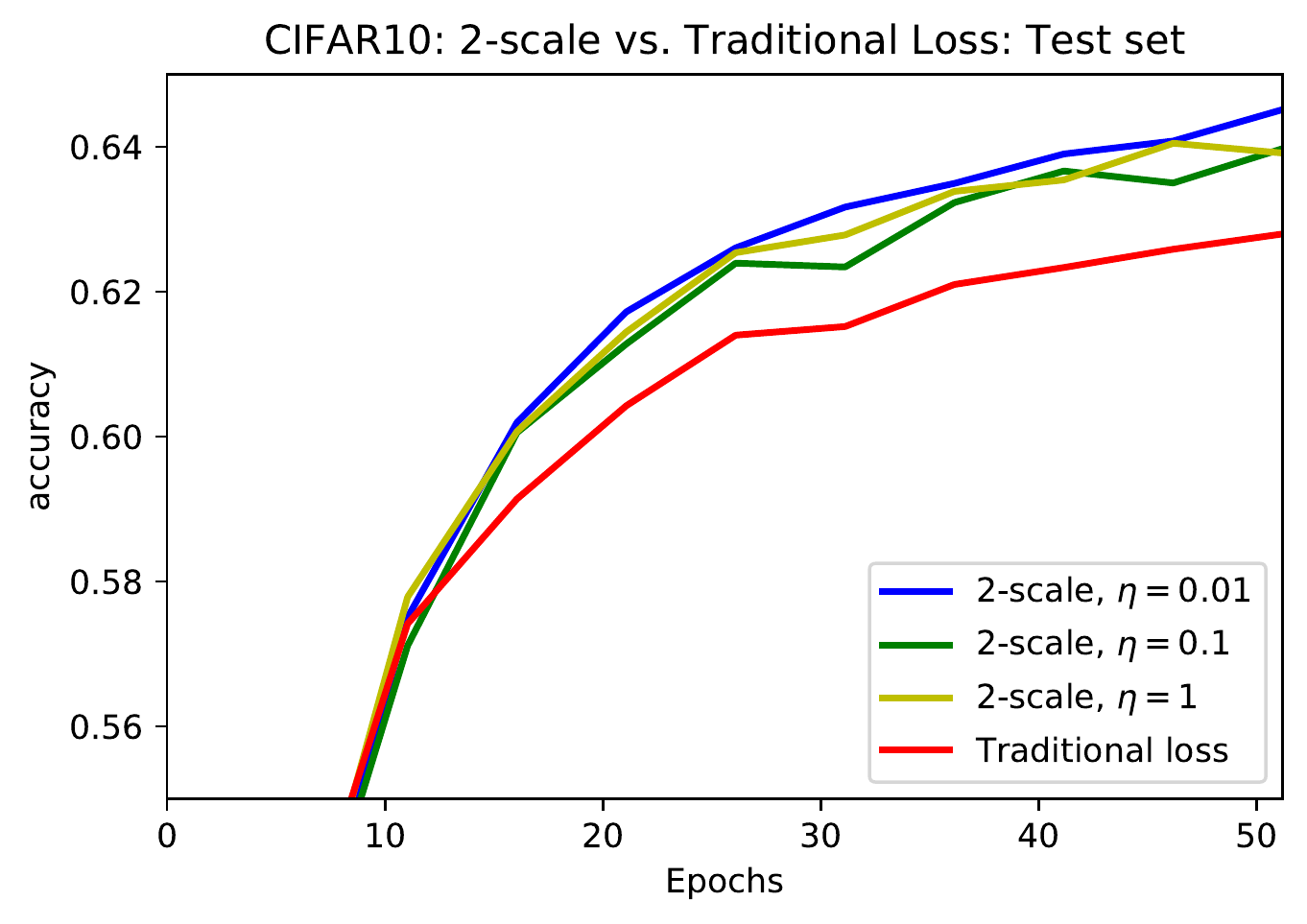}
		\caption{The accuracy throughout training on the training set (left) and the test set (right) for two-scale loss network \eqref{2_scale_loss} was compared with the traditional loss network for various values of $\eta$ on the CIFAR10 data-set.  A zoomed in version of the accuracy on the test set is also displayed (bottom). Accuracy on the test set was recorded once every 2000 iterations ($\approx$ 5 epochs) and the graphs represent the average accuracy over 10 random initial seeds. \label{Fig:2_scale_CIFAR_10_lt}}
	\end{figure}
	
	For training times of approximately 50 epochs, we observe that the two-scale loss network exhibits a large increase of 12\% in accuracy for $\eta=.01$ when compared on the training set, see Fig \ref{Fig:2_scale_CIFAR_10_lt} (left).  
	Accuracy on the test set increased by a significant amount, 1.5\%, over the same training period see Fig \ref{Fig:2_scale_CIFAR_10_lt} (right and below).  
	%
	As $\eta$ increases, the accuracy of the two-scale loss network approaches that of the traditional loss as expected.  The two-scale loss network was also tested for $\eta=10$ on CIFAR10 and the results were nearly indistinguishable compared to the traditional loss as expected.

	\begin{figure}[!t]
		\centering
		
		\includegraphics[scale=.5]{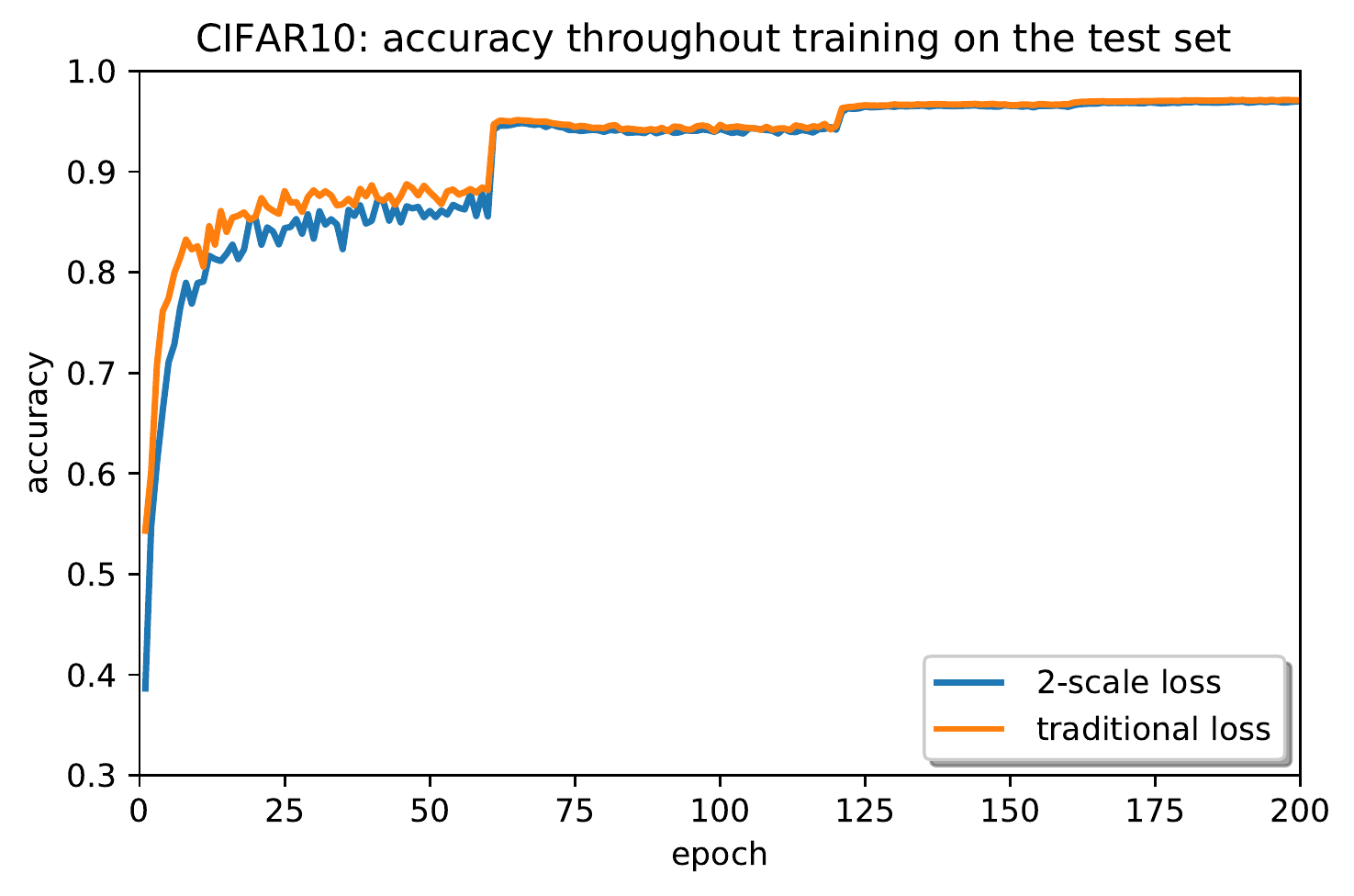}
		\includegraphics[scale=.5]{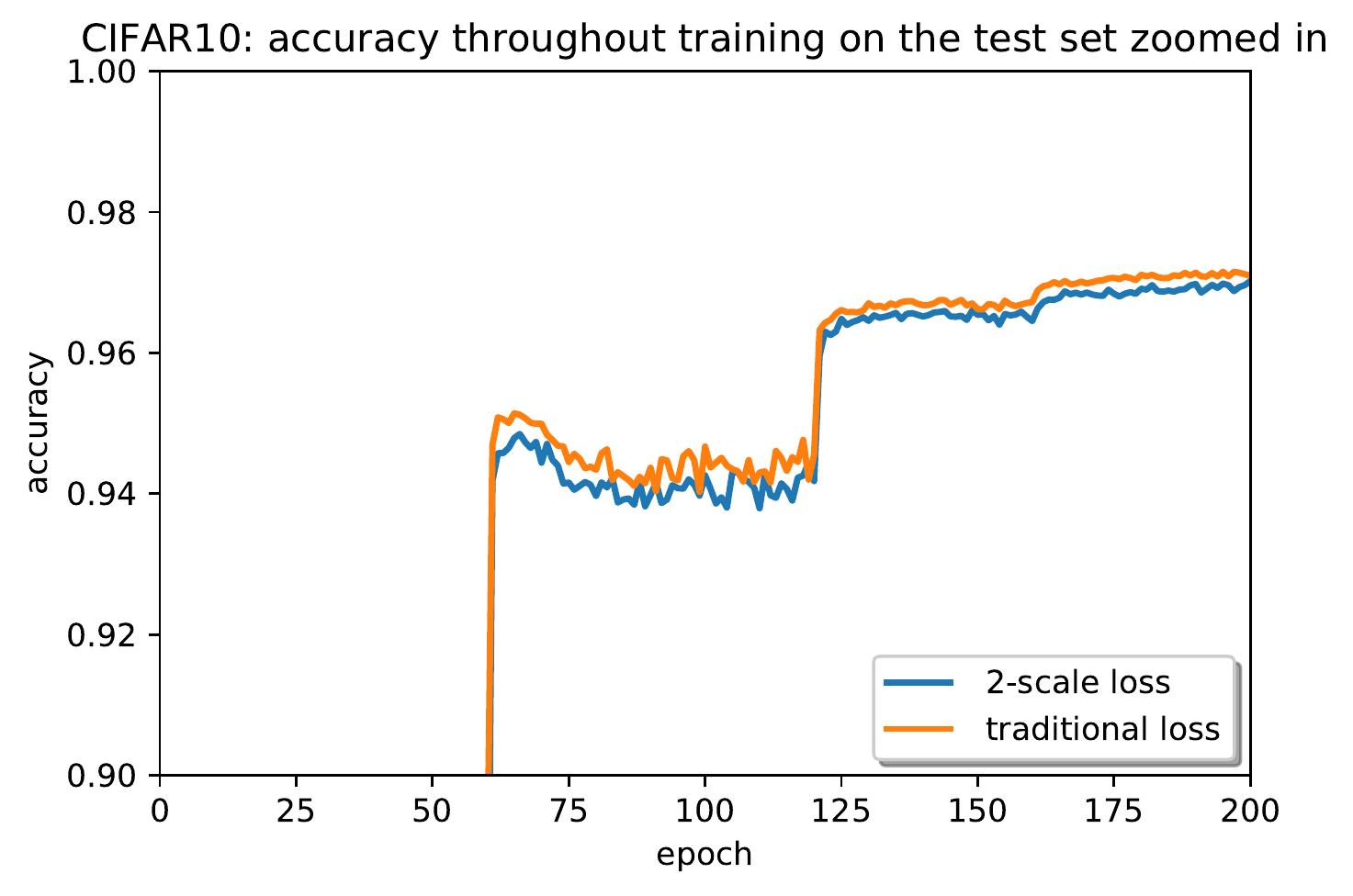}
		\caption{The accuracy throughout training on the test set for two-scale loss WRN \eqref{2_scale_loss} was compared with the traditional single scale loss WRN for various values of $\eta$ on the CIFAR10 data-set.  Accuracy on the test set was recorded approximately every 390 iterations ( 1 epoch) and the graphs represent the average accuracy of 5 random initial seeds.  The right graph is a zoomed in version of the graph on the left. \label{Fig:2_scale_CIFAR10_sam}}
	\end{figure}
\begin{figure}[!t]
		\centering
		
		\includegraphics[scale=.5]{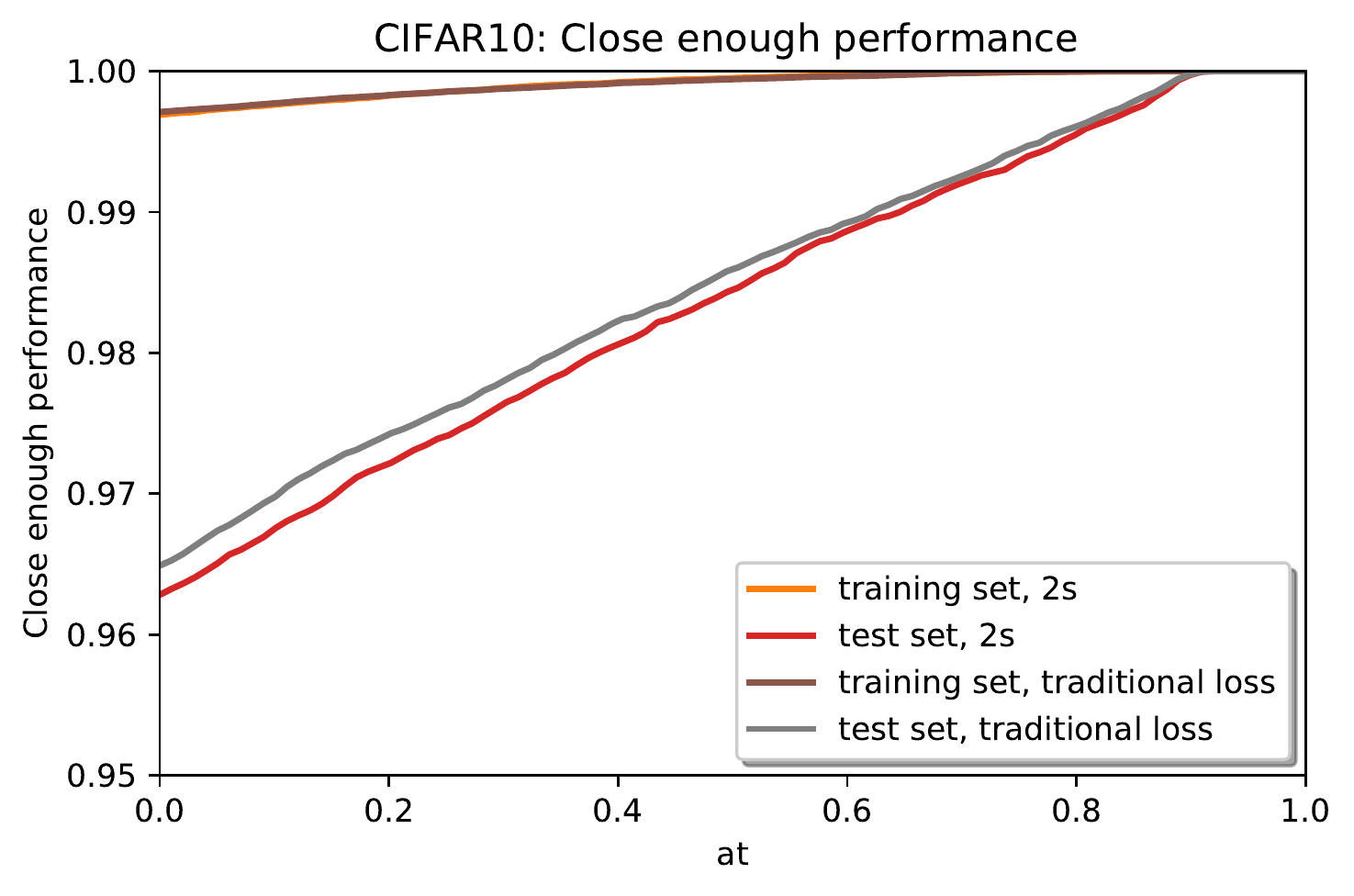}
		\caption{The close enough performance on the test set for two-scale loss WRN \eqref{2_scale_loss} was compared with the traditional single scale loss WRN for various values of $\eta$ on the CIFAR10 data-set.  The performance was recoreded at the end of training ( 200 epochs) and the graphs represent the average performance of 5 random initial seeds. \label{Fig:close_enough_CIFAR10_sam}}
	\end{figure}

	We also compare the performance for CIFAR10 using the WRN with ASAM.  We observe very high and similar accuracies on the test set for both the two-scale and single scale networks, see Fig. \ref{Fig:2_scale_CIFAR10_sam}.  The single scale network ends with a marginal increase in accuracy, $<0.2$\%. We next look at the close-enough performance to see the effect of two-scale loss on the miss-classified objects.  We observe that the gap in performance shrinks as $at$ grows but is never overcome, see Fig. \ref{Fig:close_enough_CIFAR10_sam}, contrary to what was observed in CIFAR100.  This lack of change in performance when $at$ is increased is likely due to the very high accuracy on both the test and training sets.  The high accuracy on the test set means there are less objects to capture in order for the two-scale network to improve the close-enough performance and pass the single scale network performance.  The high accuracy on the training set (>99.7\%) also shows that the two-scale network is mostly acting as a single scale network by the end of the training and will not have an significant advantage at that point.
	
	Testing was next performed on the MNIST data-set \cite{LecCor10}, a set of gray-scale, $28$x$28$ pixel images of handwritten digits.  In this case the structure of the network was two fully connected layers followed by the soft-max function \eqref{softmax}.  A total of 30 seeds were tested for approximately 2 epochs (1000 iterations of SGD).  The network was run for 30 different initial random $\alpha^0$, and the accuracy is averaged over all seeds.  
	\begin{figure}[!t]
		\centering	\includegraphics[scale=.5]{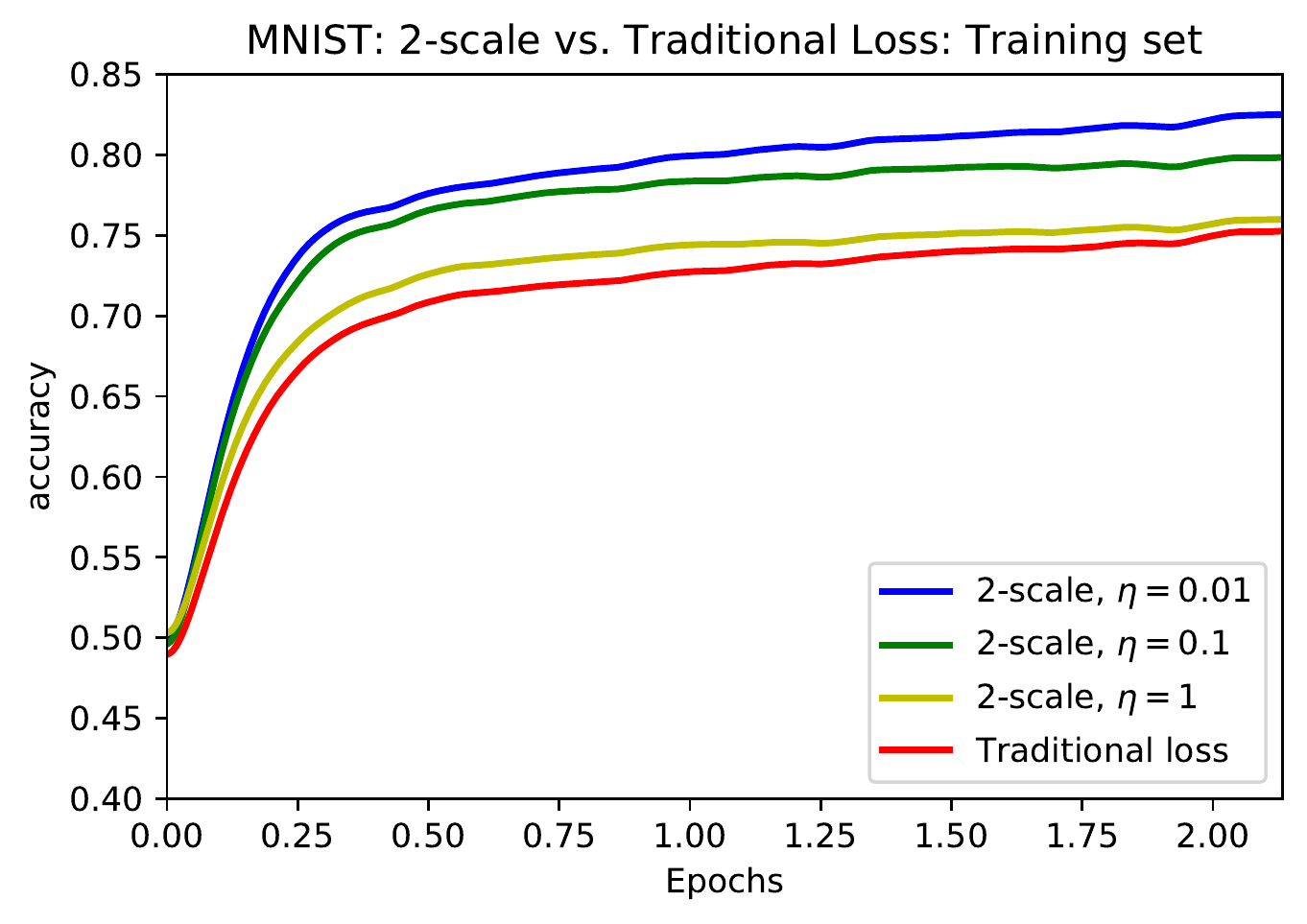}
		\includegraphics[scale=.5]{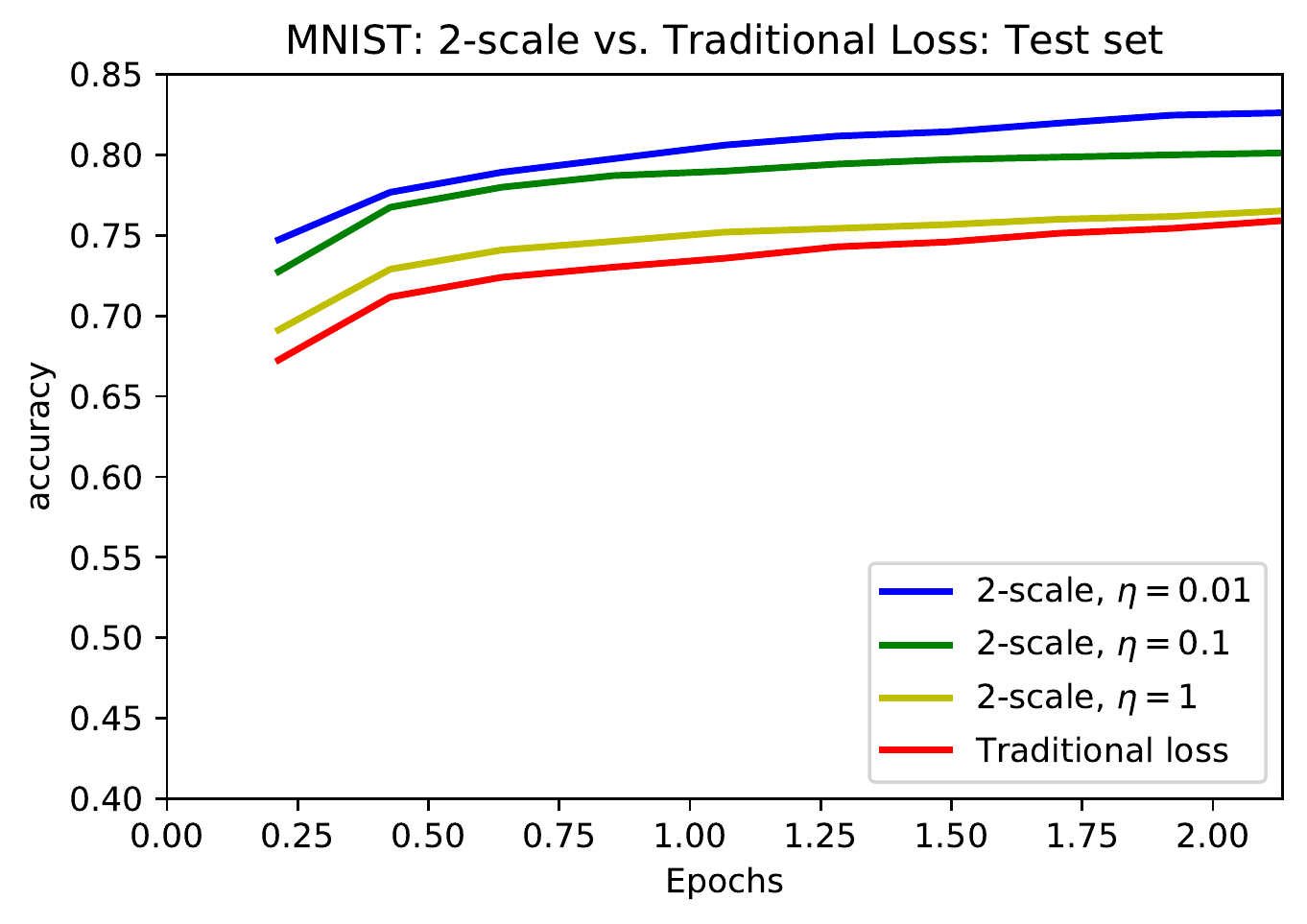}
		\caption{The accuracy throughout training on the training set (left) and the test set (right) for two-scale loss network \eqref{2_scale_loss} was compared with the traditional single scale loss network for various values of $\eta$ on the MNIST data-set.  Accuracy on the test set was recorded once every 100 iterations ($\approx$0.2 epochs) and the graph represents 30 initial seeds. \label{Fig:2_scale_MNIST_10_lt}}
	\end{figure}
	
	The two-scale loss network with $\eta=0.01$ had a larger accuracy on the training set by approximately 7.1\% after 2 epochs (1000 iterations of \eqref{GD_2s}), see Fig. \ref{Fig:2_scale_MNIST_10_lt} (left).  
	In addition the two-scale loss network had a comparable gain of overall accuracy, 6.7\%,  on the test set as well.  
	
	As observed on CIFAR10, for larger $\eta$ the behavior of the two-scale network approaches the traditional loss network.   We tested the two-scale loss network on MNIST for $\eta=10$ and $\eta=100$ as well.  For $\eta=10$ the graph of accuracy was indistinguishable from the traditional loss, as accuracy differed by less than $.1\%$ on average throughout training.  The two-scale loss network with $\eta=100$ had the exact same accuracy on each batch in training and the test set as the traditional loss for each initial $\alpha^0$ tested.  
	
	We finally note that, very early on during training (in particular the first 20 iterations), the two-scale loss networks performed worse on the training set when compared to the traditional loss on MNIST. In fact in that very early phase traditional loss outperformed two-scale loss by up to an average of 9.6\% per batch. A similar effect can be seen on CIFAR10 in the first $5$ epochs or so, so that training is  slightly slower initially for the two-scale network.
	
	To conclude the investigations in this subsection, it appears that the two-scale network also significantly outperforms the traditional loss on CIFAR10 and MNIST. For such simpler data sets, the better performance on originally poorly classified objects also translates to higher overall accuracy.
	\subsection{Numerics: Two-scale loss network with fixed scales}
	To offer a more complete analysis of possible choices of two-scale loss functions, we also investigate the so-called fixed two-scale loss function from Sec. \ref{Sec:2_scale_loss}.  We test the training algorithm defined by \eqref{GD_fixed_2s} and compare again the performance to the traditional single scale loss function.  The initial $R_0$ is determined by setting $R_0=R(\alpha^0)$ when $\alpha^0$ is chosen.  The initial choice of of parameters in \eqref{2_scale_loss} are $R_1=R_0$ and $R_2=10R_1$.
	
\begin{figure}[!t]
		\centering
		\includegraphics[scale=.5]{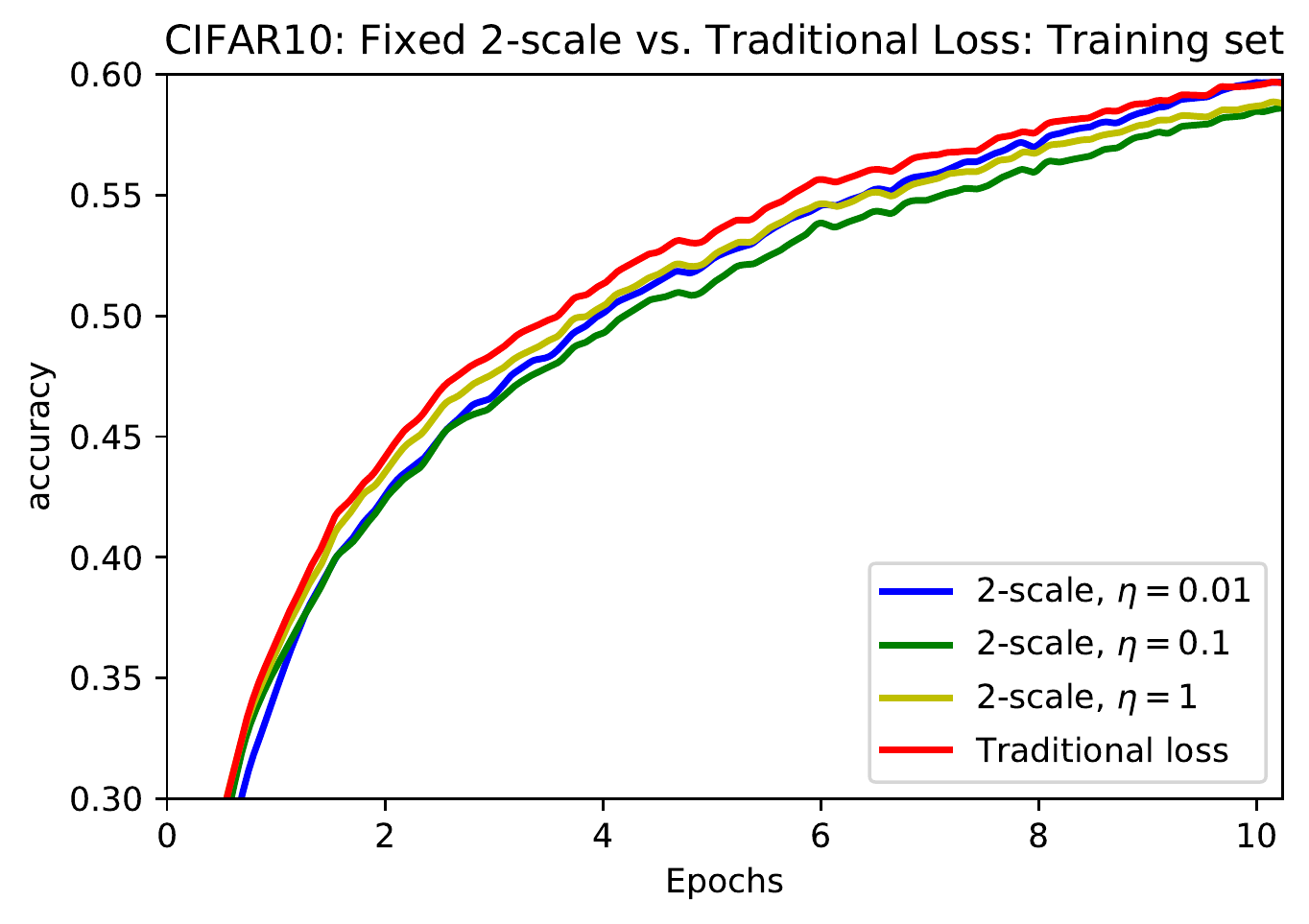}
		\includegraphics[scale=.5]{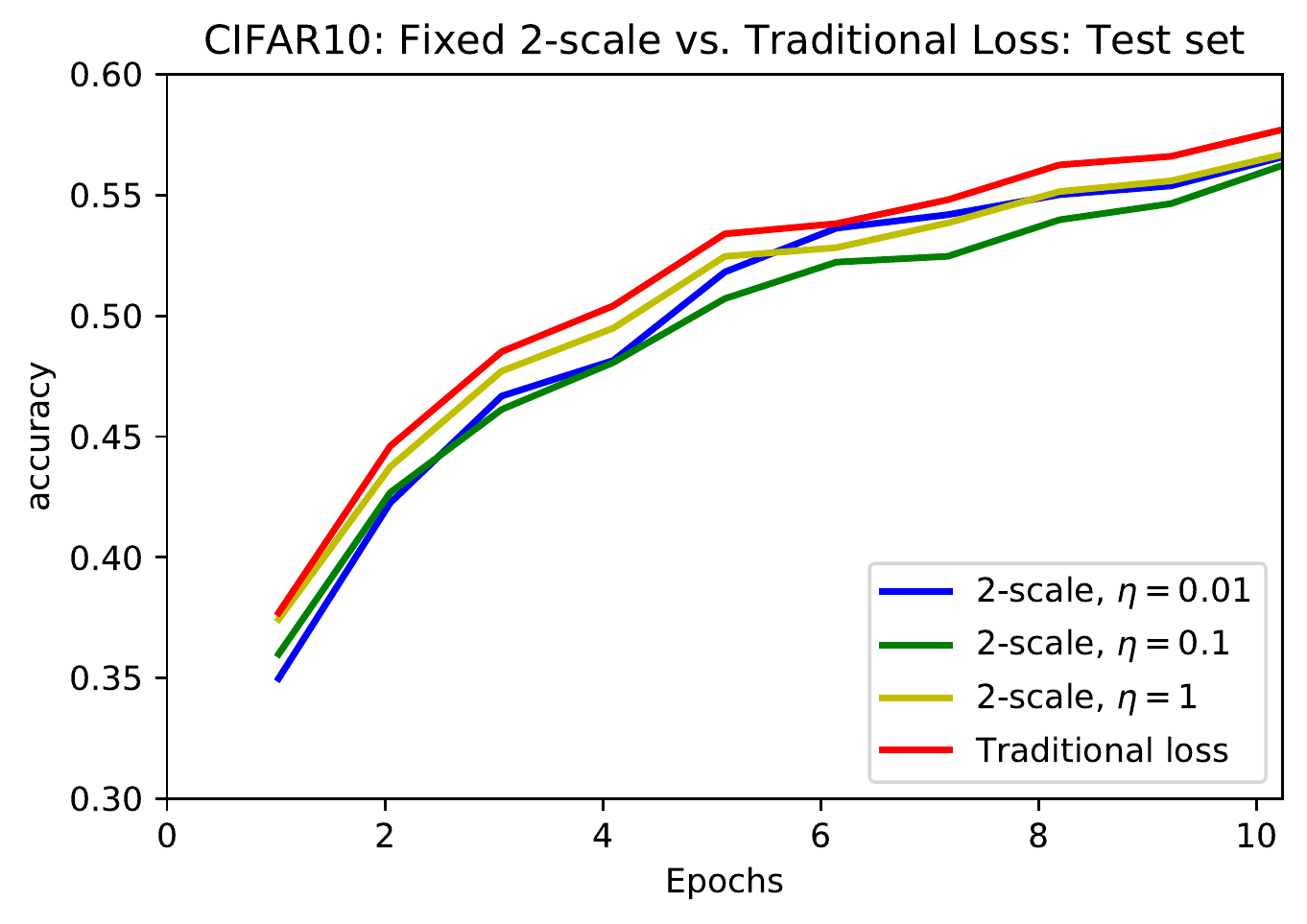}
		\caption{The accuracy on the training set (left) and the test set (right) for the fixed two-scale loss network from Sec. \ref{Sec:2_scale_loss} was compared with the traditional single scale loss function for various values of $\eta$ on the CIFAR10 data-set.  The accuracy is taken after $\approx10$ epochs of training (4000 iterations).  The graphs represent the average accuracy of 10 random initial seeds.   \label{Fig:fixed_2_scale}}
	\end{figure}

	The network is tested on CIFAR10 with the same DNN structure as in Sec. \ref{Sec:Num:2_scale}.  We observe that the the traditional loss outperforms the fixed two-scale loss network on the training set, see Fig. \ref{Fig:fixed_2_scale} (left), and even more on the testing set, see Fig. \ref{Fig:fixed_2_scale} (right).  This is of course very different from the result for the non-fixed two-scale loss network, see again \ref{Fig:2_scale_CIFAR_10_lt} (left).  Recalling that the traditional loss network outperformed the two-scale loss network, the fixed two-scale case never catches up with the performance of traditional loss and the difference between the accuracy of the networks is roughly constant.  Hence, a key takeaway is the importance of letting the two scales adapt during training.
	
	On a separate note, the behavior in terms of the choice of $\eta$ also appears more complex then before. The fixed two-scale network exhibits its worst performance for intermediary values of $\eta$ around $0.1$. As before $\eta=0.01$ leads to better results than $\eta=0.1$ but larger values of $\eta$ such as $\eta=1$ are also better, likely because they get closer to the performance of the traditional loss network. 
	For $\eta>10$ the two-scale loss network actually coincided with the gradient descent in terms of accuracy.  
	

	\subsection{Numerics: two-scale loss network with varying scale separation}
	The final loss function we consider is the two-scale separation loss network \eqref{2_scale_loss_sep}, where the traditional scale $R$ is used for the not well-classified objects and evolves via gradient descent as normal.  The scale for the well classified objects is $R_sR$, where $R_s$ denotes the separation between scales. Here $R_s$ evolves via gradient descent, using \eqref{GD_2s_sep}.  We again consider the CIFAR10 data-set and the same structure of the network as in the previous sub-sections when comparing to traditional loss.
	%
	%
	\begin{figure}[!t]
		\centering
		\includegraphics[scale=.5]{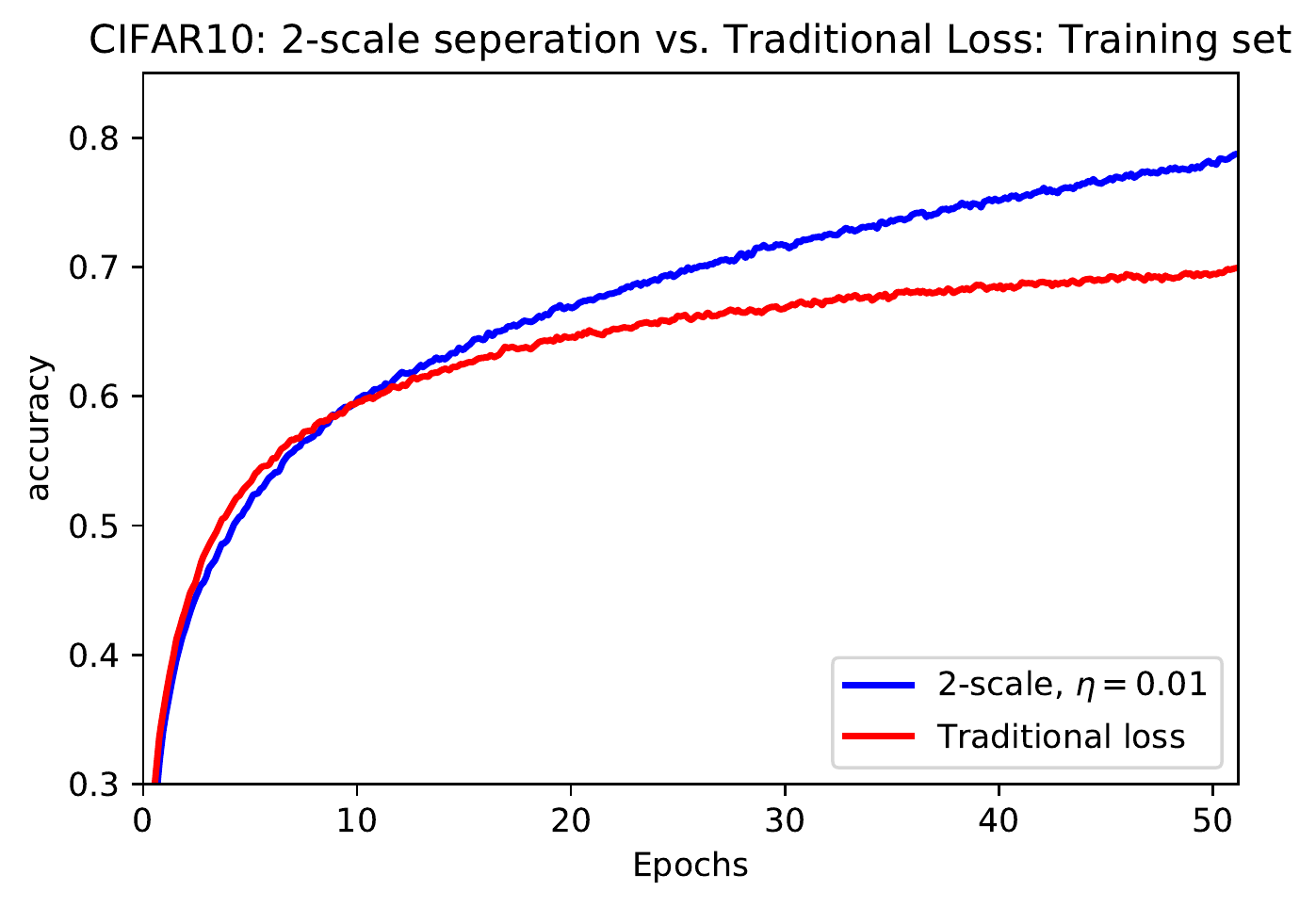}
		\includegraphics[scale=.5]{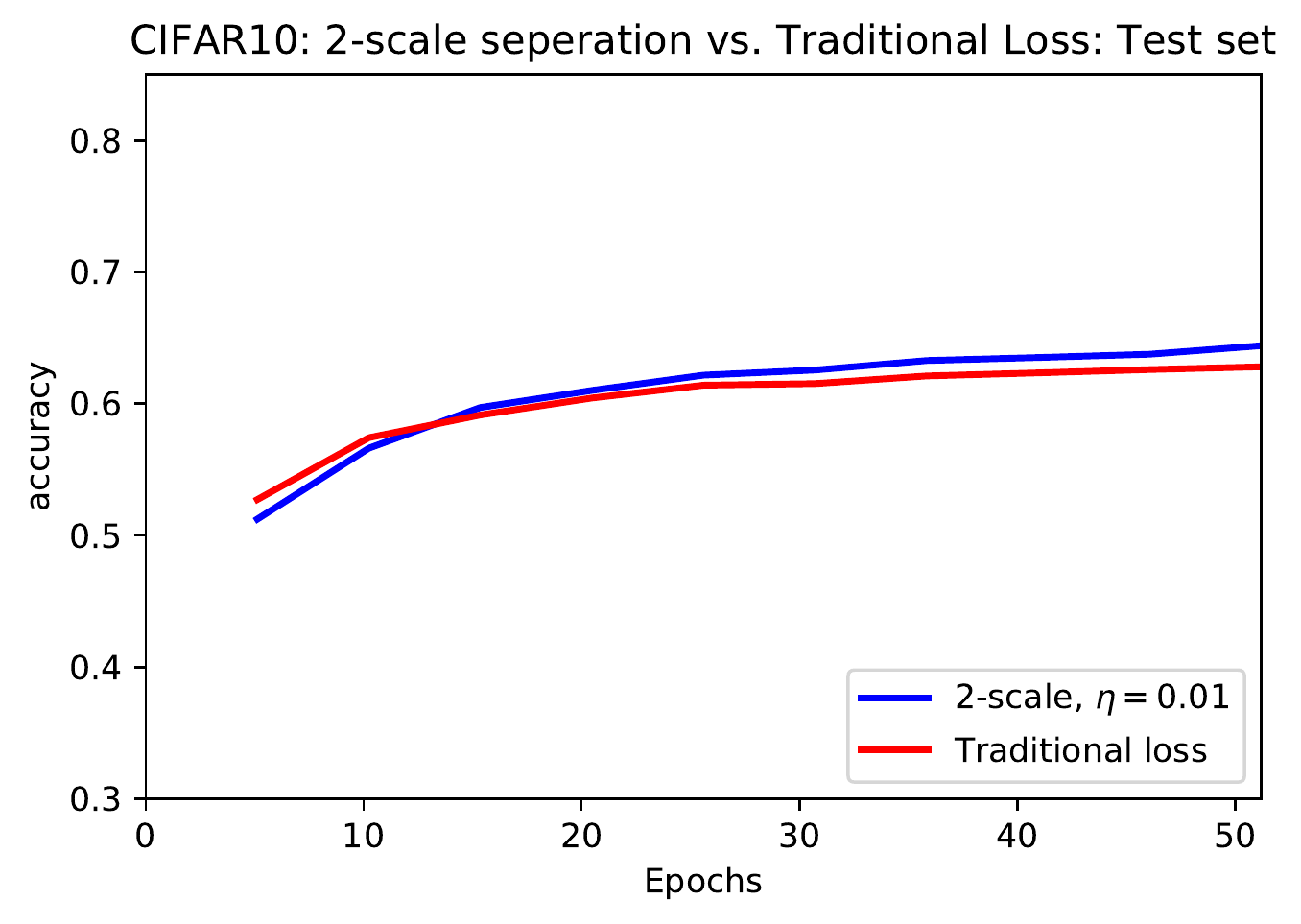}
		\caption{The accuracy on the training set (left) and the test set (right) for the two-scale loss network with varying separation from Sec. \ref{Sec:2_scale_loss} was compared with the traditional single scale loss network for various values of $\eta$ on the CIFAR10 data-set.  Accuracy on the test set was recorded once every 2000 iterations ($\approx$ 5 epochs) and the graphs represent the average accuracy over 10 random initial seeds. \label{Fig:sep_2_scale_lt}}
	\end{figure}	
	
	We observe that the two-scale separation network reaches a significantly higher accuracy on the training set and a  slightly higher accuracy on the test set, see Fig. \ref{Fig:sep_2_scale_lt}.  The behavior of the accuracy for the two-scale separation loss network is in fact very similar to that of the original two-scale loss network. This demonstrates that the main reason for the gain in accuracy is to have an adaptive ratio between the two scales.  However, the original two-scale loss network still achieved a higher accuracy than the two-scale separation loss function; approximately 1\% higher on the training set and 0.5\% higher on the test set throughout training.
	
	Overall this justifies the choice of the original two-scale loss network as the best two scale algorithms. 
	
	\section{Discussion}
	After introducing a notion of classification confidence and size of the parameters, we introduced three novel two-scale loss functions in Sec. \ref{Sec:2_scale_loss}.  The goal of these loss functions was to achieve a higher accuracy with the same DNN by focusing training on objects which are not well classified. The first two-scale loss function has independently varying scales $R_1$ and $R_2$ that were updated via SGD. The second two-scale loss function considers fixed values of $R_1$ and $R_2$.  The third two-scale loss function fixed the first scale but uses a parameter $R_s$, updated via GD, to separate the two scales \eqref{2_scale_loss_sep}.  We specifically focused on the case where the soft-max normalization and cross-entropy loss functions are used, but the two-scale loss functions can be applied to a large variety of DNN or machine learning architectures and loss functions.  
	

	After the introductions of the two-scale loss functions, numerical comparisons were made with the traditional single scale loss function on the MNIST, CIFAR10, and CIFAR100 data-sets.  We observed that the first two-scale loss network achieves a higher accuracy than the traditional loss function by a large amount on the training set on all data-sets tested. On the test sets we observed a large increase in accuracy for MNIST data-set, a smaller increase for CIFAR10, and similar accuracy on CIFAR100.
	
	We also compared the performance on CIFAR100 by using top-k performance, and an introduced measure of performance called close-enough performance.  We observed that the two-scale loss network had a significant increase of performance on both the top-n and close-enough performance measures.  The performance on the super-classes of CIFAR100 was similar to the one on the usual classes: significantly better accuracy on the training set, comparable accuracy on the testing set but with significant gains in terms of close-enough performance.
	
	The fixed two-scale loss network did not improve accuracy on CIFAR10, while the two-scale separation network did improve accuracy over the traditional loss network but did not perform as well as than the first two-scale loss network.
	
	Those results demonstrate that adaptive multi-scale loss function can provide significant increase in performance of the classical cross-entropy while only adding a minuscule amount of computational complexity and being applicable to a wide range of architectures.

	{\bf Acknowledgments} \\
	The work of L.B. and R.C. was partially supported by NSF grant DMS-2005262.  The work of P.-E. J. was partially supported by NSF grant DMS-2049020.
	
	\bibliography{bibliography}
\end{document}